\documentclass[12pt]{amsart}
\usepackage{paralist}
\def\phi{\varphi}
\def\arg{{\sf arg}}
\def\ov{\overline}

\newcommand{\C}{\mathbb{C}}
\newcommand{\CN}{\mathbb{C}^N}

\newcommand{\R}{\mathbb{R}}
\newcommand{\bR}{\mathbb{R}}
\newcommand{\N}{\mathbb{N}}
\newcommand{\Z}{\mathbb{Z}}
\newcommand{\Aut}{{\rm{Aut}_f}}

\newcommand{\al}{\alpha}
\newcommand{\bt}{\beta}
\newcommand{\gm}{\gamma}
\newcommand{\sg}{\sigma}

\newcommand{\td}{\tilde}

\newcommand{\ddop}[2]{\frac{{\rm d}#1}{{\rm d}#2}}

\newcommand{\nequiv}{{\equiv \!\!\!\!\!\!  / \,\,}}
\newcommand{\p}{\prime}
\newcommand{\tV}{\mathcal{V}}
\newlength{\extendaxesby}

\DeclareMathOperator{\id}{\sf id}

\DeclareMathOperator{\imag}{Im}
\DeclareMathOperator{\real}{Re}

\DeclareMathOperator{\ord}{ord}
\newtheorem{thm}{Theorem}
\newtheorem{lem}[thm]{Lemma}
\newtheorem{prop}[thm]{Proposition}
\newtheorem{cor}[thm]{Corollary}

\theoremstyle{definition}
\newtheorem{defin}{Definition}
\newtheorem{exple}{Example}
\newtheorem{rem}[thm]{Remark}

\def \bC{\mathbb C}

\begin{document}
\title[Degenerate real hypersurfaces in $\mathbb{C}^2$ with few
automorphisms]{
Degenerate real hypersurfaces in $\mathbb{C}^2$ with few automorphisms}%
\author{Peter Ebenfelt}%
\address{University of California, San Diego}%
\email{pebenfel@math.ucsd.edu}%
\author{Bernhard Lamel} \address{Universit\"at Wien}
\email{lamelb@member.ams.org}%
\thanks{The first author was supported in part by NSF grants DMS-0100110 and DMS-0401215. The
second author was supported by the ANACOGA network and the FWF, Projekt P17111.
The third author was supported in part by RCBS Grant of the Trinity College Dublin.
This publication has emanated from research conducted with
the financial support of Science Foundation Ireland}
\author{Dmitri Zaitsev}
\address{School of Mathematics, Trinity College, Dublin 2, Ireland}
\email{zaitsev@maths.tcd.ie}
\maketitle
\begin{abstract} We introduce new biholomorphic invariants for real-analytic hypersurfaces
  in $\bC^2$ and show how they can be used to show that a hypersurface possesses few
  automorphisms. We give
  conditions, in terms of the new invariants, guaranteeing that the stability group is
  finite, and give (sharp) bounds on the cardinality of the stability
  group in this case. We also give a sufficient condition for the
  stability group to be trivial. The main technical tool developed
  in this paper is a complete (formal) normal form for a certain class of hypersurfaces.
  As a byproduct, a complete classification, up to
  biholomorphic equivalence, of the finite type hypersurfaces in this class is
  obtained.
\end{abstract}
\section{Introduction}

Let $M$ be a germ at a point $p$ of a real-analytic hypersurface in
$\bC^2$.  An automorphism of the germ $(M,p)$ is a germ of a
biholomorphic map $H\colon (\CN,p)\to(\CN,p)$ that satisfies
$H(M)\subset M$. The set of all automorphisms of a germ $(M,p)$
forms a group under composition, called the stability group of $M$
at $p$. Endowed with the topology of uniform convergence on compact
neighbourhoods of $p$, it becomes a topological group (a sequence of
automorphisms $(H_j)$ converges to an automorphism $H$ if all of the
$H_j$ extend to a common compact neighbourhood of $p$ and converge
uniformly to $H$ on it).

It is a well known fact that a ``general'' real-analytic
hypersurface does not possess any nontrivial automorphisms. This
observation goes back to Poincar\'e \cite{Po1}. He observed that the
existence of nontrivial automorphisms imposes very strict conditions
on the coefficients of a real-analytic defining function. Indeed, if
one follows his arguments, one sees that a real-analytic hypersurface
in general position does not have any nontrivial
automorphisms! On the other hand, it is in general a hard task to
show that a given specific hypersurface has no automorphisms, as
there are no general tools available to answer this question.

In this paper, we introduce a construction that allows us, among
other things, to identify certain low order invariants associated
to a germ of a real-analytic hypersurface $(M,p)$ in $\C^2$ and to
give conditions in terms of these invariants guaranteeing that the
hypersurface does not have any nontrivial automorphisms
(Theorem~\ref{T:auttrivial}). We also give conditions guaranteeing
that the stability group is finite and provide estimates on the
number of automorphisms in this case (Theorem~\ref{T:finitedet}).
We mention here that, by means of ad hoc computations, some
explicit examples of real hypersurfaces with a finite stability
group were given in \cite{BERAJM} and \cite{St}.

The hypersurfaces we study here are Levi degenerate, that is, their
Levi forms vanish at the chosen reference point $p$. In a certain
sense, our invariants measure this vanishing in a qualitative way.
The invariants we introduce  (Theorem~\ref{T:main1}) are tensors
associated to certain lattice points $(\alpha,n,\mu) \in \N^3$,
\begin{equation}\label{E:tensor}
\mathcal L^{(\al,n,\mu)}\in{\bigotimes}^\alpha \tV_p^* \otimes
{\bigotimes}^n\bar \tV_p^*\otimes {\bigotimes}^\mu
(T^{0,1}_p\bC^2/\tV_p)^* \otimes (T^{0,1}_p\bC^2/\tV_p),
\end{equation}
where $\tV_p:=T^{0,1}_p\bC^2 \cap (\C\otimes T_pM)$. These tensors
generalize the Levi form as well as some tensors introduced by the
first author \cite{E3}.

The conditions we give in Theorems \ref{T:finitedet} and
\ref{T:auttrivial}, guaranteeing that a real hypersurface has few
automorphisms, are elementary number-theoretical relations between
the invariant lattice points $(\alpha,n,\mu)$. In order to formulate
the conditions, we have to assume that there are enough (at least
two) of these triples for a given $(M,p)$. Theorems \ref{T:finitedet} and
\ref{T:auttrivial} are given in \S\ref{S:main}, after the precise
definition of the invariants.

Our main technical result, which allows us to prove the results
described above, is a formal normal form
(Theorem~\ref{T:normalform}) for a class of hypersurfaces defined in
terms of our invariants. As in the well-known Chern-Moser normal
form, our normal form gives rise to a completely algorithmic
construction of the normalization map (by induction).

We would like to remark here that our invariants arise both in the finite type and
in the infinite type case. In the infinite type case, our results quantify the general
jet determination theorem in $\C^2$ obtained in \cite{ELZ1}; by this, we mean that
we can compute the jet order needed for the determination property
(for our class of hypersurfaces) from the invariants introduced in this paper.

Another remark in order is that in the finite type case, our
normal form also gives (in a standard way) a complete
classification,  with respect to biholomorphic equivalence, of the
hypersurfaces under consideration. Our normal form is formal, that
is, we do not prove convergence, but in the finite type case this is
not necessary, since we have the result on convergence of formal
mappings in \cite{JAMS} at our disposal. In this context, the result
in \cite{JAMS} implies that any formal invertible mapping between
two real-analytic, finite type hypersurfaces in $\bC^2$ is
convergent and, hence, yields a biholomorphism between the two
hypersurfaces. Thus, the formal classification that follows from the
normal form gives rise to a biholomorphic classification in this
setting.

Here is a short outline of this paper: in  Section \ref{S:main} we
introduce our invariants, state our main results, and discuss some
equivalent ways to define the invariants. After that, we prove the
transformation rules for our tensors in Section \ref{S:newtensors}.
Section \ref{s:normalform} is devoted to the construction of the
formal normal form, from which the main results of the paper will follow.
Finally in Section~\ref{S:Example} we show by explicit examples that our bound on the jet order
needed to determine the automorphisms is indeed achieved.

\section{Main results}\label{S:main}

Let $M\subset \bC^2$ be a germ of a real-analytic hypersurface
through $p=(p^1,p^2)\in \bC^2$. Recall that this means that $M$ is defined,
locally near $p$, by the equation $\rho=0$, where $\rho$ is a real-analytic
function near $p$ with $\rho(p)=0$ and $d\rho(p)\neq 0$. We shall
identify $\rho$ with its Taylor series
\begin{equation}\label{E:def0}
\rho(Z,\bar Z):=\sum_{I,J}\rho_{I\bar J}(Z-p)^I\ov{(Z-p)^J},
\end{equation}
where standard multi-index notation is used and the coefficients
satisfy the reality condition $\rho_{I\bar J}=\ov{\rho_{J\bar
I}}\in\bC$. More generally, we shall consider formal (not
necessarily convergent) power series (e.g.\ Taylor series of a
defining function of a germ of a smooth hypersurface) $\rho(Z,\bar
Z)$ as in \eqref{E:def0}, which satisfy the above reality condition
and the nondegeneracy condition $d\rho=\rho_{1\bar0}dZ+\rho_{0\bar
1}d\bar Z\neq 0$ at $p$, and say that $\rho$ defines a {\em formal
hypersurface} $M$ at $p\in \bC^2$.  A formal automorphism of $M$ at
$p$ is an invertible formal holomorphic mapping $(\bC^2,p)\to
(\bC^2,p)$ such that $\rho(H(Z),\overline{H(Z)})=a(Z,\bar
Z)\rho(Z,\bar Z)$ for some formal power series $a(Z,\bar Z)$ in
$Z-p$ and $\ov{Z-p}$. An invertible formal holomorphic mapping
$H\colon (\bC^2,p)\to (\bC^2,p)$ is a pair of formal holomorphic
power series $H=(H^1,H^2)$ of the form
\begin{equation}
H^j(Z)=p^j+\sum_{|I|>0}H^j_I(Z-p)^I
\end{equation}
such that
$\det \frac{\partial H}{\partial Z}(p) =
\begin{vmatrix}
  H^1_{(1,0)} & H^1_{(0,1)} \\
  H^2_{(1,0)} & H^2_{(0,1)}
\end{vmatrix}
\neq 0$. We shall denote by $\Aut(M,p)$ the group of all formal
automorphisms of $M$ at $p$.

For a real-analytic hypersurface $M$, we may choose local
holomorphic coordinates $(z,w)$ vanishing at $p$, such that $M$ is
given, locally near $p=(0,0)$, by
\begin{equation}\label{E:defeq1}
\imag w = \phi(z,\bar z, \real w)
\end{equation}
where $\phi(z,\bar z, s)$ is a (real valued) real-analytic
function in a neighborhood of $0$ in $\C^2\times\R$ satisfying
\begin{equation}\label{phi-normalization}
\phi(z, 0, s) \equiv \phi(0,\bar z, s) \equiv 0.
\end{equation}
For
a formal hypersurface $M$ through $p\in\bC^2$, the analogous
transformation, in which $\phi$ is a formal power series, is
possible by a formal holomorphic change of coordinates. Any such
coordinates $(z,w)$ (formal or local holomorphic) are called {\em
normal coordinates} for $M$ at $p$ (for more details, see e.g. \cite{BERbook}).
For each representation
\eqref{E:defeq1} of $M$, we consider the power series expansion
\begin{equation}\label{phi-expand}
\phi(z,\chi,s)=
\sum_{\al\ge 1,\,\mu\ge 0} \phi_{\al,\mu} (\chi) z^\al s^\mu.
\end{equation}

We should point out that normal coordinates $(z,w)$, as described
above, are highly non-unique. However, we will show below that for
certain lattice points $(\al,\mu)\in\mathbb N^2$, the
corresponding coefficients $\phi_{\al,\mu}(\chi)$ transform in a
particularly simple way under changes of normal coordinates (see
Proposition \ref{P:0basic}). We shall refer to such points $(\al
,\mu)$ as invariant pairs; the exact definition is given in
Definition \ref{i-pair} below.
The invariance
(that is, independence of the choice of normal coordinates) and
transformation law for the corresponding coefficient is then
established in Proposition \ref{P:0basic}
below. The lattice points which are invariant
correspond precisely to the ``lowest order coefficients'' of the
Taylor series in \eqref{phi-expand} in the sense of the following partial
ordering on $\N^2$:

\begin{equation}\begin{aligned}
  (\alpha,\mu) &\preceq (\beta,\nu) \text{ if }
  \alpha+\mu \leq   \beta +\nu \text{ and }
  \mu \leq \nu; \\
  (\alpha,\mu) &\prec (\beta,\nu) \text { if }
  (\alpha,\mu) \preceq (\beta,\nu)
  \text{ and } (\alpha,\mu)\neq(\beta,\nu).
  \label{e:precdef}\end{aligned}
\end{equation}

\begin{defin}\label{i-pair}
  A point $(\al_0,\mu_0)\in\N^2$ is called an {\em invariant pair
  associated to}
  $M$ if $\phi_{\al_0,\mu_0}(\chi)\not\equiv 0$ but
  $\phi_{\al,\mu}(\chi)\equiv 0$ for every $(\al,\mu)\prec
  (\al_0,\mu_0)$.
\end{defin}

In the following, we shall denote the set of all invariant pairs
 by $Q_{M,p}=Q_M\subset \N\times\N$.  Even though, {\it a priori}, the set $Q_M$ depends on
 the choice of normal coordinates $(z,w)$, we shall show
  (see Theorem \ref{T:main1} below) that, in fact, it does not and hence
  the set $Q_M$ is an invariant of $(M,p)$. We shall, moreover, define a
  refined invariant set $\Lambda_M \subset \mathbb N\times\mathbb
N\times\mathbb N$ as follows. For each $(\alpha,\mu)\in Q_M$, we
set
\begin{equation}
\label{E:nam} n(\alpha,\mu):=\min\left\{
n\colon
\ddop{^n \varphi_{\alpha,\mu}}{\chi^n} (0)\neq 0\right\},
\end{equation}
and define
\begin{equation}\label{E:lambdaM}
\Lambda_M:=\{(\alpha,n,\mu)\in \mathbb N\times\mathbb
N\times\mathbb N\colon (\alpha,\mu)\in Q_M\,\text{{\rm and
$n=n(\alpha,\mu)$}}\}.
\end{equation}
It is not difficult to see,
using the fact that $\phi(z,\chi,s)=\bar\phi(\chi,z,s)$, that for
any invariant pair $(\alpha,\mu)$ we have
\[ n(\alpha,\mu)\geq \alpha.\]

We are now in a position to state the main results of this paper.
Our principal {\em technical} result consists of a construction, for
each pair of points $(\alpha,n,\mu)\ne(\alpha',n',\mu')\subset
\N\times\N\times\N$ with $\alpha\neq n$, a formal normal form for
the hypersurfaces $M\subset \bC^2$ satisfying
$(\alpha,n,\mu),(\alpha',n',\mu')\in \Lambda_M$. The normal form is
described in Theorem \ref{T:normalform}. (To describe it precisely
requires distinguishing several cases and we prefer to do this at
the end of the paper.)

The normal form in Theorem \ref{T:normalform} allows us to bound the
dimension of the stability group of a real hypersurface $M\subset
\bC^2$ satisfying the condition above; it also implies a
criterion for the stability group to be trivial. Our
first result along these lines is the the following theorem, which
guarantees that the stability group can be embedded in a suitable
jet group. The construction of the
normal form also  implies, that the (formal) stability group can
be given the structure of a finite dimensional Lie group; see Theorem~\ref{T:jetparam}.
The theorem also provides bounds on the dimension of this group.

\begin{thm}\label{T:finitedet} Let $M\subset\C^2$ be a
real-analytic (or formal)
  hypersurface with  $p\in M$. Assume that
  the invariant
  set $\Lambda_M$, as defined above, contains at least two points and at least one of them, say
  $(\al,n,\mu)$, satisfies $\al \neq n$. Then,
  the group $\Aut(M,p)$ of all formal automorphisms of $(M,p)$
  embeds, via its jet evaluation, as a closed Lie subgroup of some jet group $J_p^k(\C^2)$,
which satisfies
\begin{equation}\label{E:dimauto1}
\dim_\bR\Aut(M,p)\leq1.
\end{equation}
Moreover, if either
\begin{equation}\label{E:finite1}
\alpha+n=\alpha'+n', \quad \text{for some }  (\alpha',n',\mu')\ne (\alpha,n,\mu) \in
\Lambda_M,
\end{equation}
or the number
\begin{equation}\label{E:finite2}
\frac{(\alpha'+n')\mu- (\alpha+n)\mu'}{(\alpha'+n')-(\alpha+n)}
\end{equation}
is not the same positive integer for all choices of
$(\alpha',n',\mu')\ne (\alpha,n,\mu)\in\Lambda_M$, then $\Aut(M,p)$
embeds as a finite subgroup of ${\sf U}(1)\times {\sf U}(1)$,
more precisely as a subgroup of the finite group $N_M$ described in
Remark~$\ref{R:prelimnorm}$ below. In particular,
\begin{equation}\label{E:dimauto2}
\#\Aut(M,p)\leq 2(n-\al).
\end{equation}
If $\alpha+ n \neq
\alpha' +n'$ for all $(\alpha',n',\mu')\ne (\alpha,n,\mu)\in \Lambda_M$ and the number
\eqref{E:finite2}
is the same positive integer for all $(\alpha',n',\mu')\ne (\alpha,n,\mu)\in
\Lambda_M$, then $\Aut(M,p)$ embeds into $J^k_p (\C^2)$.
\end{thm}

The given bound \eqref{E:dimauto2} on the number of automorphisms is actually sharp
as it will be demonstrated by Example \ref{ex:finauto} below.

\begin{rem} {\rm The fact that $\Aut(M,p)$ embeds into $J_p^k(\C^2)$
implies, in particular, that the automorphisms in $\Aut(M,p)$ are
determined by their $k$-jets at $p$, i.e.\ if $H,H'\in \Aut(M,p)$
and
\begin{equation}\label{E:finite3}
\partial^\alpha H(p)=\partial^\alpha  H'(p),\quad \forall
|\alpha|\leq k,
\end{equation}
then $H=H'$.}
\end{rem}

Theorem \ref{T:finitedet} is a direct consequence of the normal form
in Theorem \ref{T:normalform}. It was proved by the authors in
\cite{ELZ1} that for any real-analytic Levi non-flat hypersurface
$M\subset\bC^2$ and $p\in M$ there exists a number $k$ so that the
automorphisms in $\Aut(M,p)$ are determined by their $k$-jets at
$p$.
Moreover, if $M$ is of finite type at $p$, then $k=2$ suffices to
determine the automorphisms in $\Aut(M,p)$ (also proved in
\cite{ELZ1}).  Theorem \ref{T:finitedet} can be viewed as a
refinement of the results in \cite{ELZ1} taking into account the
finer invariants introduced in this paper.

We should also point out
that, as was shown by R. T. Kowalski \cite{Travis1}, \cite{Travis} and
the third author \cite{Zsurvey}, for any positive integer $k$, there exist
real-analytic Levi non-flat hypersurfaces $M\subset \bC^2$ with
$p\in M$ (and $M$ of infinite type at $p$) for which the
automorphisms in $\Aut(M,p)$ are uniquely determined by their
$k$-jets at $p$ but not by their $(k-1)$-jets at $p$. However, the
examples given in \cite{Travis1}, \cite{Travis} and \cite{Zsurvey}
are such that their invariant sets $\Lambda_M$
consist of a single point and hence
do not belong to the class considered in Theorem \ref{T:finitedet}.
However, in Section \ref{S:Example}, we give
another family of examples that belong to that class
(for which $\Lambda_M$ consists of at least two points) and for
which still arbitrarily high order jets are needed to determine the automorphisms.

We now come to the criterion mentioned in the introduction for $M$
to have no nontrivial automorphisms. This result is a direct
consequence of Theorem \ref{T:normalform} and the observation made
in Remark \ref{R:prelimnorm} below, and its precise formulation is
the following:

\begin{thm}\label{T:auttrivial} Let $\Lambda$ be a subset
  of $\N^3$ that contains at least two points, and one of them, say
  $(\alpha,n,\mu)$, satisfies $n \neq \alpha$.
  Assume, in addition, that
  either
\begin{equation}\label{E:finite12}
\alpha+n=\alpha'+n',
\quad \text{for some }  (\alpha',n',\mu')\ne (\alpha,n,\mu) \in \Lambda,
\end{equation}
or the number
\begin{equation}\label{E:finite22}
\frac{(\alpha'+n')\mu- (\alpha+n)\mu'}{(\alpha'+n')-(\alpha+n)}
\end{equation}
is not the same positive integer
for all choices of $ (\alpha',n',\mu')\ne (\alpha,n,\mu)\in\Lambda$.
Then all (formal or real-analytic) hypersurfaces $M$
satisfying $\Lambda_M=\Lambda$ have $\Aut(M,p) = \left\{ \id \right\}$ if
and only if
\begin{equation}\label{e:gcdcon}
  \gcd\left\{ n'-\alpha' \colon (\alpha',n',\mu') \in \Lambda \right\} = 1,
\end{equation}
there exists an even $\mu'$ with $(\alpha',n',\mu') \in \Lambda$, and
either of the following two conditions is fulfilled:
\begin{compactenum}[\rm i)]
\item $n'-\alpha'$ is even for some $(\alpha',n',\mu') \in \Lambda$
  with $\mu'$ even;
\item $n'-\alpha'$ is odd for some $(\alpha',n',\mu')\in\Lambda$ with
  $\mu$ odd.
\end{compactenum}
\end{thm}

We  note here that, in the definition of the greatest common divisor
above, we use the convention that $0$ is divisible by any integer.
Let us also note that if all the conditions of Theorem~\ref{T:auttrivial}
are fulfilled except for i) and ii), then the automorphism
group has at most $2$ elements (this follows from the explicit
form of the group $N_M$ given in \S\ref{ss:preliminarynormalization}).

We will now give some examples where Theorem~\ref{T:auttrivial}
implies the triviality of the automorphism group. In all cases, the
point $p$ is the origin.

\begin{exple} Assume that $a\geq 1$, $b\geq 1$ are integers satisfying
  $a > b +1$. Let $r\geq 0$ be an integer and assume that $p$ and $q$ are
  nonequal odd primes. Then the real hypersurface $M = M(a,b,p,q,r)$ given
  by
  \[ \imag w =(\real w)^r \left( |z|^{2a} \real{z^p} + \real w |z|^{2b} \real{z^q}
  \right)\]
  does not have any nontrivial automorphisms.
  Here $\Lambda_M=\{(a,a+p,r),(b,b+q,r+1)\}$.

  Let us check that the conditions of Theorem~\ref{T:auttrivial} are
  fulfilled. The condition
  $a> b+1$ ensures that $(a,r)$ and $(b,r+1)$ are
  invariant pairs. If $2a + p \neq 2b+q$, then the fraction in \eqref{E:finite22} is
  an integer if and only if
  \[ \frac{2a + p}{(2b + q) - (2a + p)} \in\N;\]
  but in this last fraction, the numerator is odd, while the
  denominator is even, so this is not the case.
  Since $\gcd\left\{ p,q \right\} = 1$, condition \eqref{e:gcdcon}
  is fulfilled.
  Also, either $r$ or $r+1$ is odd; so (ii) in the last condition in Theorem~\ref{T:auttrivial}
  is fulfilled.
\end{exple}

\begin{exple}
  Generalizing the last example a bit, let $a,b,r,s$ be integers,
  with $s$ being odd, satisfying $a > b+ s$, and
  $p$ and $q$ nonequal odd primes.
  Then the
  hypersurface $M = M(a,b,r,s,p,q)$ given by
  \[ \imag w =\left( \real w \right)^r \left(
  |z|^{2a} \real z^p + \left( \real w \right)^s |z|^{2b} \real z^q \right)\]
  does not have any nontrivial automorphisms.
  The invariants are given by $(\alpha,n,\mu) = (a,a+p,r)$ and
  $(\alpha',n',\mu') = (b, b+q ,r+s)$ in this example.
  Again, let us check that the conditions of Theorem~\ref{T:auttrivial} are
  fulfilled.
  The conditions $a > b+s$ ensures that $(a,r)$ and $(b,s)$ are invariant pairs.
  The fraction in  \eqref{E:finite22} is an integer if and only if
  \[ \frac{s(2a+p)}{(2b+q) - (2a + p) } \in \N;\]
  just as in the preceding example, this is never the case.
  The reasoning of the preceding example also applies to the verification of
  the last two conditions of Theorem~\ref{T:auttrivial}.
\end{exple}


\begin{exple}
  Let $a > 3$. The hypersurface $M$ given by
  \[ \imag w = (\real w) |z|^2 (\real z) + |z|^{2a} \]
  has no nontrivial automorphisms.

  Let us again check the conditions of Theorem~\ref{T:auttrivial}.
  The invariants are given by $(\alpha,n,\mu) = (1,2,1)$ and $(\alpha',n',\mu')= (a,a,0)$.
  The fraction condition~\eqref{E:finite22}
  is satisfied, since $1< \frac{2a}{2a-3} < 2$. Furthermore, the condition for
  the $\gcd$ is satisfied since $n -\alpha =1$, which is odd, so  ii) holds. Actually,
  i) holds also, since $0$ is even.
\end{exple}

\begin{rem}
  Let us remark that the automorphism groups remain trivial even if we allow
  higher order terms to appear in the Taylor expansion of the hypersurfaces in
  the
  examples above. Here, higher order terms has to be understood in the sense of
  the partial ordering used in Definition~\ref{i-pair}; that is, the terms of the form
  $z^\beta \bar z^n (\real w)^\nu$ with $ (\alpha,\mu) \prec (\beta,\nu)$
and $(\alpha',\mu') \prec (\beta,\nu)$.
\end{rem}

We now give an example of a hypersurface having the maximal number
of automorphisms allowed by the bound in \eqref{E:dimauto2}.

\begin{exple}
  \label{ex:finauto} Let $a,k,q$ be positive integers such that $q$
  is odd and
  $2a+k$ is not divisible by 3.
   Consider the hypersurface
  given by
\[ \imag w = (\real w)^{q}\real{z}^k \left( |z|^{2a}(\real w)^{2} +
 |z|^{2(a+3)} \right). \]
  The invariants are $(\alpha,n,\mu)=(a,a+k,q+2)$ and
 $(\alpha',n',\mu')=(a+3,a+3+k,q)$.
 The fraction \eqref{E:finite2} is given by
 \[ \frac{(\alpha' + n')\mu - (\alpha+n) \mu'}{(\alpha' + n') - (\alpha+n)} = q+2 + \frac{2a+k}{3},\]
 which is not an integer since we assume that $3$ does not divide $2a+k$.
Hence, by Theorem \ref{T:finitedet},
 the group $\Aut(M,0)$ is finite and satisfies the estimate
 \eqref{E:dimauto2}. Since the biholomorphisms
 \[ H_{\ell,+} (z,w): = \left( e^{\frac{2 \ell \pi i }{k}} z, w  \right),
 \quad H_{\ell,-}(z,w): =  \left( e^{\frac{ 2 \ell \pi i }{k}} z, - w  \right), \quad 0
 \leq\ell \leq k-1.
 \]
 are all in $\Aut(M,0)$, we conclude that we actually have equality
 in \eqref{E:dimauto2} in this case; also, this implies that
 we have $\Aut(M,0) = \left\{ H_{\ell,-}, H_{\ell,+} \colon 0\leq \ell \leq k-1 \right\}$.

\end{exple}

\begin{exple}
  It is natural to ask what kind of role the fraction condition plays. An example showing that a
  one-parameter family of higher jet parameters of
  the order predicted by Theorem~\ref{T:finitedet} can be really needed if $2\le
  \frac{(\alpha' + n')\mu - (\alpha+n) \mu'}{(\alpha' + n') - (\alpha +n)} \in \N$
is given in \S\ref{S:Example}
  below; here, we give an example showing  that there is a dependence on a one-parameter family
  of first order jets in case
 $\frac{(\alpha' + n')\mu - (\alpha+n) \mu'}{(\alpha' + n') - (\alpha +n)} =1$,
  thus establishing that the dimension bound \eqref{E:dimauto1} is sharp also in that case.

  Let $\alpha\geq 1$ and $\mu \geq 2$ be integers, and let $p$ be a positive
  integer such that $\alpha - \mu + 1 \geq p$. Then the hypersurface $M$ given by
  \[ \imag w = (\real w)^{\mu} \left( |z|^{2\alpha} + (\real w)^{\mu -1}
  \real( z^{\alpha-\mu+1-p} \bar z^{3\alpha+\mu - 1 +p})\right) \]
  has the one-parameter family of biholomorphisms given by
  \[ H_r (z,w) = \left( r^{1-\mu} z, r^{2\alpha} w \right), \quad
   0 \neq r \in \R.\]
  Actually, we note here that the normal form given in Theorem~\ref{T:normalform} and the
  observations made in \S\ref{ss:preliminarynormalization} imply that
  the automorphism group is generated by $H_r$ and the
  discrete group of rotations
  \[ \left( z,w \right) \mapsto \left( e^{\frac{\pi i j}{n+p+\mu -1}}z,w \right) , \quad
  0\leq j < 2(n+p+\mu- 1).\]
\end{exple}

We will now discuss Definition~\ref{i-pair} in some detail.
Let us first note that for fixed normal
coordinates, we could also decompose $\phi (z, \chi, s) =
\sum_{\bt\ge
  1, \, \mu\ge 0} \psi_{\bt,\nu} (z) \chi^\bt s^\nu$. Let us show that
this representation leads to the same invariant pairs, where we use
the analogous definition as above.  Since $\phi$ is real valued, $\phi
(z,\chi,s) = \bar \phi (\chi,z,s)$.  Hence,
\[ \sum_{\al,\mu} \phi_{\al,\mu} (z)\chi^\al s^\mu = \phi (z,\chi,s) =
\bar \phi (\chi,z,s) = \sum_{\bt,\nu} \bar \psi_{\bt,\nu} (z) \chi^\bt
s^\nu, \] and so
\begin{equation}\label{E:reality1}
\phi_{\al,\mu}(z) = \bar \psi_{\al,\mu}(z) \text{ for all } \al,
  \,\mu.\end{equation}
We shall also use a different kind of defining equation for $M$
for which the calculations turn out to be simpler. It is well
known (the reader can consult for example the book by Baouendi,
Ebenfelt and Rothschild \cite{BERbook}
for details) that $M$ can also be given in the form
\begin{equation}\label{E:normality}
w = Q(z,\bar z,\bar w),
\end{equation}
where $Q(z,\chi,\tau)$ is a holomorphic function (or a
formal power series) in a neighborhood
of $0$ in $\C^3$ satisfying
\begin{equation}\label{normality}
Q(z,0,\tau) \equiv Q(0,\chi,\tau) \equiv \tau,
\quad Q(z,\chi,\bar Q(\chi,z,\tau)) \equiv \tau.
\end{equation}

We may also use the function $Q$ to define invariant pairs as
follows. We decompose $Q(z,\chi,\tau)$ as
\begin{equation}\label{E:Qisequal0}
  Q(z,\chi,\tau) = \tau + \sum_{\al\ge 1,\, \mu\ge 0} q_{\al,\mu}
  (\chi) z^\al \tau^\mu
 = \tau + \sum_{\bt\ge 1, \, \mu\ge 0} r_{\bt,\nu} ( z) \chi^\bt
 \tau^\nu,
\end{equation}
and define invariant pairs to be the minimal
coefficients in the ordering used in Definition \ref{i-pair},
with $\phi_{\alpha,\mu}$ replaced by either $q_{\al,\mu}(\chi)$ or
$r_{\al,\mu}(\chi)$. The
equivalence of all these definitions is a consequence of the
following lemma.

\begin{lem}\label{equivalence}
  In given normal coordinates $(z,w)$, let $M$ be defined by each of the
  equations \eqref{E:defeq1} and \eqref{E:normality}. Then for every
  pair $(\al_0,\mu_0)\in\N^2$, the following properties are
  equivalent:
\begin{enumerate}
  \item[{\rm(i)}] $\phi_{\al,\mu}(\chi)\equiv 0$ for $(\alpha,\mu)\prec(\alpha_0,\mu_0)$;
  \item[{\rm(ii)}] $q_{\al,\mu}(\chi)\equiv 0$  for $(\alpha,\mu)\prec(\alpha_0,\mu_0)$;
  \item[{\rm(iii)}] $\psi_{\al,\mu} (z)\equiv 0$ for $(\alpha,\mu)\prec(\alpha_0,\mu_0)$;
  \item[{\rm(iv)}] $r_{\bt,\nu}(z)\equiv 0$  for $(\alpha,\mu)\prec(\alpha_0,\mu_0)$.
\end{enumerate}
If {\rm (i)}--{\rm (iv)} are satisfied, then
\begin{equation}\label{q-phi-r}
q_{\al_0,\mu_0}(\chi)\equiv 2i \phi_{\al_0,\mu_0}(\chi) \equiv
-\bar r_{\al_0,\mu_0}(\chi).
\end{equation}
\end{lem}

\begin{proof}
  The proof is based on the identity
\begin{equation}\label{q-phi}
\frac{Q(z,\chi,\tau)-\tau}{2i} \equiv \phi\Big(z,\chi, \frac{Q(z,\chi,\tau)+\tau}{2}\Big)
\end{equation}
which is a consequence of \eqref{E:defeq1} and \eqref{E:normality}.
Let $(\al_0,\mu_0)$ satisfy (ii) and assume that
$\phi_{\al,\mu}\not\equiv 0$ for some $(\al,\mu)$ as in (i).
Then in the expansion of the left-hand
side of \eqref{q-phi} the coefficient of $z^\al \tau^\mu$ (as a
function of $\chi$) is zero whereas,
in view of (ii) and \eqref{phi-normalization}, only the term $\tau$ can contribute
in the first expansion \eqref{E:Qisequal0} of $Q(z,\chi,\tau)$ on the right-hand side.
Hence the corresponding coefficient on the right-hand
side is $\phi_{\al,\mu}$ which is not zero, a contradiction.  Similar
computation of the coefficients of $z^{\al_0} \tau^{\mu_0}$ yields
the first identity in \eqref{q-phi-r}.

Conversely, suppose that (i) holds but $q_{\al,\mu}(\chi)\not\equiv 0$
for some $(\al,\mu)$ as in (ii). Then such a pair
$(\al,\mu)=(\al_1,\mu_1)$ can be chosen such that (ii) holds with
$(\al_0,\mu_0)$ replaced by $(\al_1,\mu_1)$.  The above argument shows
that the first identity in \eqref{q-phi-r} holds with $(\al_0,\mu_0)$ replaced by
$(\al_1,\mu_1)$ which contradicts the assumption (i). Hence (ii) also
holds as required.

We have already shown (see \eqref{E:reality1}) that (i) and (iii) are
equivalent. The proof that (iii) and (iv) are equivalent is exactly
the same as above, using the expansions in terms of $\psi$ and $r$.
\end{proof}

For each $(\alpha,n,\mu)\in
\Lambda_M$ with the notation above, we define a tensor
\begin{equation}\label{E:tensor1}
\mathcal L^{(\al,n,\mu)}\in{\bigotimes}^\alpha \tV_p^* \otimes
{\bigotimes}^n\bar \tV_p^*\otimes {\bigotimes}^\mu
(T^{0,1}_p\bC^2/\tV_p)^* \otimes (T^{0,1}_p\bC^2/\tV_p),
\end{equation}
where $\tV_p$ denotes the space of $(0,1)$ vectors at
$p$ which are tangent to $M$,
in local coordinates as follows. First, observe that $\tV_{p}$
is spanned, in normal coordinates
$(z,w)$ for $M$ at $p$, by $\partial/\partial \bar z$ and the
normal space $T^{0,1}_p\bC^2/\tV_p$ is spanned by the projection
of $\partial/\partial\bar w$. By choosing $\partial/\partial \bar
z$ as a basis for $\tV_{p}$ and $\partial/\partial \bar w$ mod
$\partial/\partial \bar z$ as a basis for $T^{0,1}_p\bC^2/\tV_p$,
we identify these two spaces with $\bC$ and define
\begin{equation}\label{E:deftensor}
\mathcal
L^{(\al,n,\mu)}(a_1,\ldots,a_\alpha,b_1,\ldots,b_n,c_1,\ldots,c_\mu):=
\frac{1}{n!}\ddop{^n q_{\alpha,\mu}}{\chi^n} (0) a_1\ldots a_\alpha b_1\ldots
b_n c_1\ldots c_\mu.
\end{equation}
The following result yields a preliminary
classification of real hypersurfaces $M\subset \bC^2$ in terms of
the set $\Lambda_M$ and also proves the fact that \eqref{E:deftensor} indeed
defines a tensor, as claimed in \eqref{E:tensor1}.

\begin{thm}\label{T:main1} Let $M\subset \bC^2$ be a real-analytic
(or formal)
  hypersurface, $p\in M$, and $(z,w)$ normal coordinates for $M$ at
  $p$. Then the set $\Lambda_M$ defined by
  \eqref{E:lambdaM} and, for each $(\al,n,\mu)\in \Lambda_M$,  the tensor
  $\mathcal L^{(\al,n,\mu)}$ defined by
  \eqref{E:deftensor} are independent of the choice of normal
  coordinates. That is, if $M^\prime\subset \bC^2$ is another
  real-analytic (or formal) hypersurface given in normal coordinates $(z^\prime,
  w^\prime)$ at $p'\in M'$ by $w^\prime = Q^\prime (z^\prime, \bar
  z^\prime, \bar w^\prime)$ and $(z',w')=(F,G)$ is a (formal)
  biholomorphic map
  sending $M$
  into $M'$ with $(F(0),G(0))=(0,0)$, then $\Lambda_M=\Lambda_{M'}$
  and
  for each
  $(\alpha,n,\mu)\in \Lambda_M$
\begin{equation}\label{E:transtens}
\ddop{^n q_{\alpha,\mu}}{\chi^n} (0) \overline{G_w}(0)=
\ddop{^n {q'}_{\alpha,\mu}}{\chi^n} (0) F_z(0)^\alpha
\overline{F_z}(0)^n\overline{G_w}(0)^\mu.
\end{equation}
\end{thm}

\begin{rem} We should point out that any CR diffeomorphism between
two $C^\infty$-smooth real hypersurfaces induces a formal
invertible map between the corresponding formal hypersurfaces (see
e.g. \cite{BERbook}). Hence, the set $\Lambda_M$ and the
corresponding tensors \eqref{E:deftensor} are in fact CR
invariants.
\end{rem}

Theorem \ref{T:main1} will be proved in section \ref{S:newtensors}. We
shall in fact show the more general result that for $(\alpha,\mu)\in
Q_M$, the whole coefficient $q_{\alpha,\mu}(\chi)$ transforms like a
family of tensors (at least for $(\alpha,\mu)\neq (1,0)$); see
Proposition \ref{P:0basic}
for the precise statement.

Before proceeding further, we have two remarks: First, as
is easily seen, the set $\Lambda_M$ (or, equivalently, the set
$Q_M$) always contains at least one point unless $M$ is Levi flat,
i.e. $Q\equiv 0$. Indeed, if $q_{\beta,\nu}\not\equiv 0$, for some
$(\beta,\nu)$, then there is at least one $(\alpha_1,\mu_1)\in
Q_M$ such that $(\alpha_1,\mu_1)\prec (\beta,\nu)$.
 Second, we also point out that if
$(1,1,0)\in\Lambda_M$, then $\Lambda_M$ does not contain any other
point and the corresponding tensor $\mathcal L^{(1,1,0)}$ is
simply the Levi form; thus
our tensors generalize the Levi form
for Levi degenerate hypersurfaces.

\section{Invariant tensors and their transformation}\label{S:newtensors}
\subsection{Transformation rule and invariance of the tensors}\label{ss:definition1} Let
$M\subset\C^2$ be formal hypersurface (e.g.\ coming from a germ of
a real-analytic one) and $p\in M$. We shall keep the notation
introduced in section \ref{S:main}. Thus, $(z,w)$ will be normal
coordinates for $M$ at $p=(0,0)$ and we shall assume that $M$ is
defined by the equation \eqref{E:defeq1} or in complex form by
\eqref{E:normality}. We decompose $Q(z,\chi,\tau)$ in
two ways as in \eqref{E:Qisequal0}.

Let us recall that $M$ is of {\em finite type} at $p=(0,0)$ if and only if
there exists a $\gm$ such that $q_{\gm, 0 }(\chi) \nequiv 0$.
In this case we
will denote the smallest such $\gm$ by $\gm_0$. It is easy to
check that $\gm_0$ is independent of the choice of normal
coordinates; cf.\ also Proposition \ref{P:0basic} below.
If $M$ is of infinite type, we set $\gm_0 = \infty$.

If $H=(F,G)$ is a formal invertible map taking $M$ into another
hypersurface $M^\prime$ which is given in normal coordinates
$(z^\prime, w^\prime)$ by $w^\prime = Q^\prime (z^\prime, \bar
z^\prime, \bar w^\prime)$, then (by definition)
\begin{equation} \label{E:HmapsMintoMprime}
 G(z,w) = Q^\prime (F(z,w),  \bar F(\chi, \tau), \bar G(\chi,
 \tau)),
\end{equation}
when $w=Q(z,\chi,\tau)$. By the normality condition
\eqref{E:normality}, setting $\chi = 0$, $w = \tau$, we have that
\begin{equation}\label{E:Gisequal}
  G(z,w) = Q^\prime ( F(z,w), \bar F(0,w), \bar G (0,w)).
\end{equation}
Note that this implies that $G(z,0) \equiv 0$.
Putting $w =Q(z,\chi,\tau)$ in the right-hand sides of
\eqref{E:HmapsMintoMprime} and \eqref{E:Gisequal}
and equating them and using \eqref{E:Qisequal0} we obtain
\begin{multline}\label{E:start}
  \bar G\big(0,Q(z,\chi,\tau)\big) + \sum_{\al\ge 1, \mu\ge 0} q_{\al,
    \mu}^\prime \big(\bar F (0,Q(z,\chi,\tau)\big) \big(
  F(z,Q(z,\chi,\tau))\big)^\al \big( \bar G(0,Q(z,\chi,\tau))\big)^\mu
  \\ = \bar G (\chi,\tau) + \sum_{\al\ge1, \mu\ge0} q_{\al, \mu}^\prime
  \big(\bar F (\chi,\tau)\big) \big(F(z,Q(z,\chi,\tau))\big)^\al
  \big(\bar G(\chi,\tau)\big)^\mu,
\end{multline}
where the $q_{\al, \mu}^\prime$ are defined as in \eqref{E:Qisequal0},
with $Q^\prime$ replaced by $Q$.

Now, recall from the previous section (see Definition \ref{i-pair}
and the paragraph following it) the definition of the set
$Q_M\subset \N\times \N$ as the set of all invariant pairs
 of $M$. The following proposition proves
the invariance of $Q_M$ (hence justifying the terminology "invariant
pairs") and shows how the $q_{\al,\mu}(\chi)$, for $(\al,\mu)\in
Q_M$, transform under formal changes of normal coordinates.

\begin{prop}\label{P:0basic}
  Let $M,M'\subset\bC^{2}$ be formal hypersurfaces, each given
  in normal coordinates at $0\in M$ and $0\in M'$ respectively, and
  $H=(F,G)\colon (\bC^2,0)\to (\bC^2,0)$ a formal
 invertible map sending $M$ into $M'$.  Assume that $(\al_0,\mu_0)$ is an
  invariant pair for $M$. Then it is also an invariant
  pair for $M'$ and
\begin{equation}\label{E:0basicone}
\begin{aligned}
  q_{\al_0,\mu_0}(\chi)= &\, q'_{\al_0,\mu_0}\big(\bar F(\chi,0)\big)
  F_z(0)^{\al_0} \ov{G_{w}}(0)^{\mu_0-1},\
  \text{\rm for } (\al_0,\mu_0)\ne (1,0)\\
  q_{1,0}(\chi) =&\, q'_{1,0}\big(\bar F(\chi,0)\big)\big(
  F_z(0)+F_w(0)q_{1,0}(\chi)\big)\ov{G_{w}}(0)^{-1},\ \text{{\rm
      for $(\alpha_0,\mu_0)=(1,0)$}}.
\end{aligned}
\end{equation}
\end{prop}

\begin{proof}
We begin by assuming that
\begin{equation}\label{E:assu}
q'_{\al,\mu}\equiv 0 \quad \text{ for every }(\al,\mu)\prec (\al_0,\mu_0).
  \end{equation}
 We expand both sides of \eqref{E:start} into  a Taylor series in $z$ and
  $\tau$ and identify the coefficients of $z^{\al_0}\tau^{\mu_0}$,
  using the definition of invariant pairs and Lemma~\ref{equivalence}.
  The first term on the left-hand side of \eqref{E:start} has the Taylor expansion
\begin{equation}\label{E:Gtaylorexpanded}
  \bar G \big(0,Q(z,\chi,\tau)\big) = \sum_{l=1}^\infty \frac{1}{l!}
  \ov{G_{w^l}}(0) \Big(\tau + \sum_{\al, \mu} q_{\al,\mu} (\chi) z^\al   \tau^\mu\Big)^l.
\end{equation}
The general form of a term in the expansion of $\left(\tau + \sum_{\al, \mu}
  q_{\al,\mu} (\chi) z^\alpha \tau^\mu\right)^l$ is, up to a binomial
factor, either
\begin{equation}\label{E:gentermGlhs}
q_{\gm_1,\sg_1}(\chi) \cdots q_{\gm_{l_1},\sg_{l_1}} (\chi)
z^{\gm_1 +\dots + \gm_{l_1}} \tau^{l-l_1  + \sg_1 +\dots +
\sg_{l_1}},
\end{equation}
for some $1\leq l_1 \leq l$, or $\tau^l$.
Since $(\al_0,\mu_0)$ is invariant,
if the term \eqref{E:gentermGlhs} is not 0,
then for each $1\leq j\leq l_1$ either
$\gamma_j+\sigma_j\geq\alpha_0+\mu_0$ or $\sigma_j\geq\mu_0$. We
conclude that
\begin{equation}\label{E:gequiv}
\bar G(0,Q(z,\chi,\tau))\sim \bar G(0,\tau)+\ov{G_w}
(0)q_{\alpha_0,\mu_0}(\chi)z^{\alpha_0}\tau^{\mu_0},
\end{equation}
where $\sim$ means equality modulo terms of the form
$z^\gamma\tau^\sigma$ with either $\gamma+\sigma>\alpha_0+\mu_0$
or $\sigma>\mu_0$. Thus, the only term of the form
$z^{\alpha_0}\tau^{\mu_0}$
in the expansion \eqref{E:Gtaylorexpanded}
is $\ov{G_w}
(0)q_{\alpha_0,\mu_0}(\chi)z^{\alpha_0}\tau^{\mu_0}$.

Let us examine the sum on the left in \eqref{E:start}. An argument
similar to the one above shows that
\begin{equation}\label{E:fequiv}
F(z,Q(z,\chi,\tau))\sim F(z,\tau)+
F_w(0)q_{\alpha_0,\mu_0}(\chi)z^{\alpha_0}\tau^{\mu_0},
\end{equation}
where $\sim$ has the same meaning as in \eqref{E:gequiv}. Using
also \eqref{E:assu} and that $q'_{\alpha,\mu}(0)=0$
(by normality of the coordinates $(z,w)$)
for every $(\al,\mu)$, it is not difficult to see that
the expansion of
$$q_{\al, \mu}^\prime (\bar F (0,Q(z,\chi,\tau)) \Bigl(
F(z,Q(z,\chi,\tau))\Bigr)^\al \Bigl( \bar
G(0,Q(z,\chi,\tau))\Bigr)^\mu$$
cannot contribute a term with $z^{\alpha_0}\tau^{\mu_0}$.

Let us now examine the right-hand side of \eqref{E:start}. The
first term $\bar G(\chi,\tau)$ cannot contribute a term with
$z^{\alpha_0}\tau^{\mu_0}$ since $\alpha_0\geq 1$. In the sum on
the right of \eqref{E:start}, terms of the form
$$q_{\al, \mu}^\prime (\bar F (\chi,\tau) \Bigl(
F(z,Q(z,\chi,\tau))\Bigr)^\al \Bigl( \bar
G(\chi,\tau)\Bigr)^\mu,$$ with $\mu>\mu_0$, cannot contribute,
since as already noted above, $\bar G(\chi,0) \equiv0$. We conclude that
the contribution from the sum on the right can only come from the
terms involved in the expansion of
\begin{equation}\label{right-exp}
q'_{\al_0,\mu_0} \big(\bar F(\chi,0)\big) \big(F_z(0)z+
F_w(0)(\tau+q_{\al_0,\mu_0}(\chi)z^{\al_0}\tau^{\mu_0})\big)^{\al_0}
\big(\ov{ G_w}(\chi,0)\tau\big)^{\mu_0}
\end{equation}
and hence equals
$$q'_{\al_0,\mu_0} (\bar F(\chi,0)) (F_z(0)z)^{\al_0}(\ov{
G_w}(\chi,0)\tau)^{\mu_0} ,$$ if $(\al_0,\mu_0)\neq (1,0)$ and
$$q'_{1,0} (\bar F(\chi,0)) (F_z(0)
+ F_w(0) q_{1,0}(\chi))z,$$  if $(\al_0,\mu_0)=(1,0)$. Thus,
\eqref{E:0basicone} follows, provided that we show that $\ov{ G_w}
(\chi,0) \equiv \ov{G_w} (0,0)$ holds when $\mu_0\ge 1$. To see
this, we set $z=0$ in \eqref{E:start} and obtain, in view of
\eqref{normality},
\begin{multline}\label{E:startz0-0}
  \bar G(0,\tau) + \sum_{\al, \mu} q_{\al, \mu}^\prime \big(\bar F
  (0,\tau)\big) \big( F(0,\tau)\big)^\al \big( \bar G(0,\tau)\big)^\mu
  \\ = \bar G (\chi,\tau) + \sum_{\al, \mu} q_{\al, \mu}^\prime
  \big(\bar F (\chi,\tau)\big) \big(F(0,\tau)\big)^\al \big(\bar
  G(\chi,\tau)\big)^\mu.
\end{multline}
The desired property now follows by differentiating
\eqref{E:startz0-0} with respect to $\tau$, setting $\tau=0$, and
using the fact that $q'_{1,0}\equiv 0$ if $\mu_0\ge 1$.  The proof
of \eqref{E:0basicone} is complete and, hence also the proof of
Proposition \ref{P:0basic}, under the assumption \eqref{E:assu}.

To complete the proof, we must show that \eqref{E:assu} holds. If
\eqref{E:assu} would not hold, there would exist an invariant pair
  $(\al_0',\mu_0')\ne(\al_0,\mu_0)$ for $M'$ with $\al_0'+\mu_0' \le \al_0+\mu_0$ and
  $\mu_0'\le \mu_0$ (i.e.\ with $(\al_0',\mu_0')\prec(\al_0,\mu_0)$).
Now the corresponding assumption \eqref{E:assu} holds with $M$ and
$M'$ exchanged and the above proof yields the identities
\eqref{E:0basicone} with $H$ replaced by $H^{-1}$ and
$(\al_0,\mu_0)$ by $(\al_0',\mu_0')$. In particular, these
identities imply that $q_{\al_0',\mu_0'}\not\equiv 0$, which
contradicts the fact that $(\al_0,\mu_0)$ is an invariant pair for
$M$. Hence the assumption \eqref{E:assu} must hold and the proof is
complete.
\end{proof}

 We
are now in a position to prove Theorem \ref{T:main1}.

\begin{proof}[Proof of Theorem $\ref{T:main1}$] The invariance of the set $Q_M$
 follows directly from Proposition \ref{P:0basic}. Next, if $(\alpha,\mu)$
 is an invariant pair for $M$, then $n=n(\alpha,\mu)$ is also
 invariant. Indeed, if $H=(F,G)$ is a formal change of normal
 coordinates as in Proposition \ref{P:0basic}, then $\chi\mapsto\bar F(\chi,0)$ is a
 formal change of variables near
$\chi=0$ (since $\bar G(\chi,0)=0$) and thus it follows from
\eqref{E:0basicone}
that the order of vanishing at 0 of $q_{\alpha_0,\mu_0}(\chi)$ and
that of $q'_{\alpha_0,\mu_0}(\chi')$ are the same.
 The transformation rule
  \eqref{E:transtens} follows by Taylor expanding \eqref{E:0basicone}
   in $\chi$. \end{proof}

\subsection{An identity for the transversal component of a mapping}
\label{ss:mappings1} We shall keep the notation introduced above.
Thus, $M\subset \C^2$ and $M^\prime\subset \C^2$ denote two formal
hypersurfaces with points $p\in M$ and $p\in M^\prime$, and
$H$ denotes a formal invertible map at $p$ with $H(p)=
p^\prime$ and $H( M)\subset M^\prime$. We choose normal
coordinates $(z,w)$ for $M$ vanishing at $p$ and $(z^\prime,
w^\prime)$ for $M^\prime$ vanishing at $p^\prime$ and write
$H(z,w) = (F(z,w),G(z,w))$.

Our goal in this section is to derive a certain identity for the
transversal component $G$ of a formal change of coordinates $H$ as
above, which ensures that the normality of coordinates is preserved
by the change of coordinates induced by $H$. This is done by an
application of the invariants introduced above. We formulate this as
a lemma:
\begin{lem}\label{L:zwnormalifandonlyif}
  Let $M^\p$ be a formal hypersurface in $\C^2$, given in normal
  coordinates $(z',w')$ vanishing at $p'_0\in M'$, and let
  $H=(F,G)\colon (\bC^2,p)\to (\bC^2,p'_0)$
   be a formal invertible map.
  Let $M$ be the formal
  hypersurface $H^{-1} (M^\p)$ and $(z,w)$ be local coordinates
 vanishing at $p\in M$ so that $z^\p =
  F(z,w)$, $w^\p = G(z,w)$. Then there exist universal polynomials
  $p_k$ such that the coordinates $(z,w)$ are normal for $M$ at $p$ if and
  only if, for every $k$,
\begin{equation}\label{E:zwnormalifandonlyif}
  G_{w^k} (z,0) - \ov{ G_{w^k}} (0) =
  p_k \Bigl( \ov{ G_{w^j} }(0), F_{w^j}(0), F_{w^l} (z,0), \left(
    {r^\p_{\bt,\nu}}^{(s)} (F(z,0))\right) \Bigr),
\end{equation}
where ${r^\p_{\bt,\nu}}^{(s)}$ are the terms in the expansion
analogous to \eqref{E:Qisequal0} of the defining equation for $M'$
and  $j\leq k-m_0+1$, $l\leq k - m_0$, $\bt +\nu +s \leq k$, and
$m_0$ is defined by
\[ m_0 = \min \left\{ \alpha+\mu \colon
q_{\alpha,\mu}^{\prime} \not\equiv 0 \right\}\ge 1.\]
\end{lem}
\begin{proof} It is well known that $(z,w)$ are normal coordinates
  for $M$ at $p$ if and only if $\rho (z,0,w,w) = 0$ for some - and hence any
  - (possibly complex and formal) defining function $\rho(z,\bar z, w,\bar w)$ for $M$ (see
  e.g. \cite{BERbook}, Prop. 4.2.3). In our case, a defining function
  for $M$ is given by
$$\rho(z, \bar z,w,\bar w) := G(z,w) - Q^\p
  (F(z,w), \bar F (\bar z,\bar w) ,\bar G( \bar z,\bar w)).$$
Using the second expansion in \eqref{E:Qisequal0}
we see
  that $(z,w)$ are normal coordinates for $M$ if and only if
  \begin{equation}\label{E:Gisequal2}
  G(z,w) = \bar G (0,w) + \sum_{\bt,\nu} r_{\bt,\nu}^\p (F(z,w)) \bar
  F(0,w)^\bt \bar G(0,w)^\nu.
\end{equation}
Note that the  integer $m_0$  defined above is invariant,
and can be equivalently defined by
\begin{equation}
  \label{E:monedef'}
  m_0 = \min\{ \al + \mu : (\al,\mu)\in Q_M^\prime\}
= \min\{ \al + \mu : q_{\al,\mu}^\prime\not\equiv 0 \}
= \min\{ \bt + \nu : r_{\bt,\nu}^\prime\not\equiv 0 \}.
\end{equation}
The minimum power of $w$
appearing in the sum on the right-hand side of
\eqref{E:Gisequal2} is $m$,
hence we obtain
\begin{equation}\label{E:Gwisconstant'}
  G_{w^k} (z,0) = \ov{ G_{w^k} }(0), \quad k< m_0.
\end{equation}
To obtain a formula for $G_{w^k} (z,0)$ for $k\geq m_0$ we expand
$r_{\bt,\nu}^\p (F(z,w))$ in $w$:
\begin{equation}\label{E:r(F)expanded'}
  r_{\bt,\nu}^\p (F(z,w)) =
\sum_{n=0}^{\infty}
        a_n\big((F_{w^l} (z,0))_{l\le n}, \big(r'^{(l)}_{\bt,\nu}
        (F(z,0))\big)_{l\leq n}
        \big)
        w^n,
\end{equation}
where the $a_n$ are universal polynomials---that is, they do not
depend on either $H$ or $M^\prime$. Substituting
this in \eqref{E:Gisequal2} we obtain \eqref{E:zwnormalifandonlyif}.
\end{proof}

\section{A normal form for certain hypersurfaces in
  $\C^2$}\label{s:normalform}
\subsection{Preliminary normalization}\label{ss:preliminarynormalization}
Our goal in this section is to make a preliminary normalization of
the defining equation \eqref{E:normality} of $M$ and work out the
restrictions on the first order jets of those mappings $H=(F,G)$
that respect this normalization. First, we shall think of $M'$ as
being given,  fix an invariant pair $(\al,\mu)\neq (1,0)$ of $M'$ (i.e.\ an
element of $Q_{M'}$) assuming it exists, and find a biholomorphic mapping $H=(F,G)$,
preserving normal coordinates, such that $M=H^{-1}(M')$ is given
by \eqref{E:normality} with $q_{\al,\mu} (\chi) = 2i\lambda
\chi^{n}$ (which corresponds to
$\phi_{\al,\mu}(\chi)=\lambda\chi^n$ by \eqref{q-phi-r}), where $n
= n(\al,\mu)$ and $|\lambda|=1$. Indeed, by the construction of $n$, we
can write $q^\prime_{\al,\mu} (\chi^\prime ) = ic (\psi
(\chi^\prime))^{n}$, where $\psi(\chi^\prime ) = \chi^\prime +
O(({\chi^\prime})^2)$, and $c\in \bC$.  The map
$H(z,w):=\big(\bar\psi^{-1}(F_z(0)z),G_w(0)w\big)$, with
$F_z(0)\in\bC^*:=\bC\setminus\{0\}$ and
$G_w(0)\in\bR^*:=\R\setminus\{0\}$ considered as free parameters, preserves normality of the
coordinates and, in view of \eqref{E:0basicone}, yields
\begin{equation}\label{E:normalalpha1}
  q_{\al,\mu} (\chi) =
ic |F_z (0)|^{2\al} \overline{F_z} (0)^{n -\al} G_w  (0)^{\mu
-1}\chi^{n}.
\end{equation}

We will from now on assume that the chosen invariant pair $(\alpha,\mu)$ satisfies
$n - \alpha \neq 0$. Then choosing $F_z(0)\in\bC^*:=\bC\setminus\{0\}$ and
$G_w(0)\in\bR^*:=\R\setminus\{0\}$ suitably, we obtain
$q_{\al,\mu} (\chi)=2i \chi^n$.
For the transformations $H=(F,G)$ respecting
the normalization $q_{\al,\mu} (\chi) = 2i  \chi^{n}$, i.e.\
sending $M$ with $q_{\al,\mu} (\chi) = 2i \chi^{n}$ into
$M'$ with $q'_{\al,\mu} (\chi') = 2i (\chi')^{n}$, where
$(\al,\mu)$ is an invariant pair for both $M$ and $M'$, the
identity \eqref{E:0basicone} implies that $\bar F (\chi,0)$ is linear, i.e.\
\begin{equation}\label{2-trans}
\bar F (\chi,0) =
  \overline{F_z}(0)\chi.
\end{equation}
Now substituting \eqref{2-trans} into \eqref{E:0basicone} yields
\begin{equation}\label{E:normalalpha1.2}
 1 = |F_z(0)|^{2\al} \ov{ F_z }(0)^{n - \al} \ov{ G_w}  (0)^{\mu-1}.
\end{equation}
For integers $\alpha,n,\mu$ we shall denote by
$C(\alpha,n,\mu)$ the set of parameters $(F_z(0),G_w(0))$ satisfying \eqref{E:normalalpha1.2},
that is,
\begin{equation}\label{e:defC}
C(\alpha,n,\mu) = \left\{
(\lambda,r) \in \C^* \times \R^* \colon
\lambda^\alpha {\bar \lambda}^{n} r^{\mu-1} =1\right\}.
\end{equation}
From \eqref{E:normalalpha1.2} we obtain restrictions on $F_z(0)$
and $G_w(0)$ depending on the integers $\al,\mu,n$.  For the
arguments we obtain
\begin{equation}\label{arg}
(n-\al)\;\arg(F_z(0)) + (\mu-1)\;\arg(G_w(0)) = 0 \mod 2\pi,
\end{equation}
where $\arg(G_w(0))$, in view of \eqref{E:Gwisconstant'}, is either $0$ or
$\pi$. Moreover, the absolute value of $F_z(0)$ is related to that
of $G_w(0)$ by
\begin{equation}\label{abs}
|F_z(0)|=|G_w(0)|^{\frac{1-\mu}{\al+n}}.
\end{equation}

We now let $({\tilde\alpha},{\tilde\mu})$ denote any other invariant pair, and
work out a normalization of
the coefficient of $\chi^{\tilde n}$ in $q_{{\tilde\alpha},{\tilde\mu}}(\chi)$.
Observe that if there are two invariant pairs, then $(1,0)$ cannot
be one (i.e.\ $q_{1,0}\equiv 0$). Thus, the transformation rule
\eqref{E:0basicone} for the pair $({\tilde\alpha},{\tilde\mu})$ is necessarily
\[ q_{{\tilde\alpha},{\tilde\mu}} (\chi) = q^\p_{{\tilde\alpha},{\tilde\mu}}
(\ov{ F}(\chi,0)) F_z (0)^{{\tilde\alpha}} \ov{ G_w}(0)^{{\tilde\mu - 1 }}.
\]
Denote the coefficient of $\chi^{ {\tilde n}}$ (where $ {\tilde n}=
n({\tilde\alpha},{\tilde\mu})$) on the left hand side by $ce^{i\theta}$ with
$c >0$, and the
coefficient of $(\chi^\p)^{{\tilde n}}$ on the right hand side by
$c^\p e^{i \theta'}$. Comparing the coefficient of $\chi^{{\tilde n}}$ on both sides,
we obtain $ c e^{i\theta}  = c^\p e^{i \theta'} \ov{ F_z}(0)^{{\tilde n}}
F_z(0)^{{\tilde\alpha}} \ov{ G_w} (0)^{{\tilde\mu - 1}}$. Taking absolute values,
and using \eqref{abs},
we obtain
\[ c = c^\p |F_z(0)|^{{\tilde\alpha} + {\tilde n}} |G_w (0)|^{{\tilde\mu} -1} =
c^\p |G_w (0)|^
{\frac{(1-\mu)({\tilde\alpha} + {\tilde n}) - (1 -{\tilde\mu})(\al + n)}{\al + n}}.
\]

We now distinguish two cases:
\begin{equation}\begin{aligned} \label{e:frac}
  \text{Case A:}\quad{\tilde{\mu}\left( \alpha + n \right) - \mu \left(
\tilde{\alpha} + \tilde{n}\right)} &= {\left( \alpha+ n \right) -
\left( \tilde{\alpha} + \tilde{n} \right)}, \\
  \text{Case B:}\quad{\tilde{\mu}\left( \alpha + n \right) - \mu \left(
\tilde{\alpha} + \tilde{n}\right)} &\neq {\left( \alpha+ n \right) -
\left( \tilde{\alpha} + \tilde{n} \right)} .
\end{aligned}
\end{equation}

In Case A, the norm $c$ of the coefficient of $\chi^{\tilde{n}}$ is
an invariant of $M$. In Case B, we can change the norm of this
coefficient and normalize it by requiring that $c = 2$.
Note that in this case, the normalization fixes $|G_w(0)|$, and thus
by \eqref{abs} also $|F_z (0)|$. Thus, we now consider formal
invertible mappings $H=(F,G)$ respecting the normalization
$q_{\al,\mu}(\chi)=2i\chi^n$ and
$q_{{\tilde\alpha},{\tilde\mu}}(\chi)= i c \varepsilon \chi^{{\tilde
n}}+O(\chi^{{\tilde n}+1})$, with $|\varepsilon|=1$, where $c$ in
Case A is the real positive invariant introduced above or $2$ in
case B.

The possible values of the unimodular number
$\varepsilon$ in the normalization of $q_{{\tilde\alpha},{\tilde\mu}}$ are
restricted to a discrete subset of the unit circle ${\sf U}(1)\subset\C$, namely
\begin{equation}\label{E:normeps}
  \varepsilon=\varepsilon_0 \frac{\bar \gamma^{{\tilde n}-{\tilde\alpha}}}
  {|\gamma|^{{\tilde n}-{\tilde\alpha}}}
  \frac{\delta^{\tilde \mu - 1}}{|\delta|^{\tilde \mu -1}},
\end{equation}
where $\varepsilon_0$ is any particular unimodular number such
that $q_{{\tilde\alpha},{\tilde\mu}}$ can be normalized as
$i c \varepsilon_0\chi^{{\tilde n}}+O(\chi^{{\tilde n}+1})$ and $(\gamma,\delta)$
range over the set $C(\alpha,n,\mu)$
given by \eqref{e:defC}. We shall further restrict the
pairs $(F_z(0),G_w(0))$ by normalizing $\varepsilon$ in
\eqref{E:normeps} so that $\arg{\varepsilon}\in [0,2\pi)$ is as
small as possible. This choice of $\varepsilon$ is an invariant of
$M$ and the pairs $(F_z(0),G_w(0))$ which preserve the
normalization
\begin{equation}\label{E:qamnorm}
q_{\al,\mu}(\chi)=2i\chi^n,\quad
q_{{\tilde\alpha},{\tilde\mu}}(\chi)=ic\varepsilon\chi^{{\tilde n}}+O(\chi^{{\tilde n}+1}),
\end{equation}
where ${\tilde n}:=n({\tilde\alpha},{\tilde\mu})$, $\varepsilon \in
{\sf U}(1)$ is the invariant just described, and $c$ is the invariant
introduced above in Case A or $c=2$ in Case B, are precisely those,
which belong to $C(\alpha,n,\mu)\cap
C(\tilde{\alpha},\tilde{n},\tilde{\mu})$. We summarize the
above normalization in the following proposition:

\begin{prop}\label{P:prelimnorm} Let $M\subset \bC^2$ be a
formal hypersurface with $p=0\in M$. Assume that there are two
invariant pairs $(\al,\mu)\neq ({\tilde\alpha},{\tilde\mu})$ and that
$n:=n(\al,\mu)\neq \al$. Then, there are normal coordinates
$(z,w)\in\bC^2$ for $M$ at $p$ such that $M$ is given there by
\eqref{E:normality} with $q_{\al,\mu}$ and $q_{{\tilde\alpha},{\tilde\mu}}$
satisfying \eqref{E:qamnorm}. Here $\varepsilon\in{\sf U}(1)$ is the
invariant defined above and $c$ is the invariant introduced above
in Case A in \eqref{e:frac} or is normalized
to be $2$ in Case B.
Moreover, if $M$ is normalized by
\eqref{E:qamnorm}, then any
formal invertible mapping $H=(F,G)\colon (\bC^2,0)\to
(\bC^2,0)$ sending $M$ into $M'$, where $M'$ is also normalized by
\eqref{E:qamnorm}, satisfies $(F_z(0),G_w(0))\in C(\alpha,n,\mu)
\cap C(\tilde{\alpha}, \tilde{n} ,\tilde{\mu})$
 and $F(z,0)=F_z(0)z$. Here $C$ is defined by \eqref{e:defC}.
\end{prop}

\begin{rem}
  We remark here that the conditions \eqref{E:finite1}--\eqref{E:finite2} and
  \eqref{E:finite12}--\eqref{E:finite22} in
  Theorems~\ref{T:finitedet} and \ref{T:auttrivial} respectively guarantee
  that with the assumptions of those Theorems, we are always in Case B.
\end{rem}
\begin{rem}\label{R:prelimnorm}
In this Remark, we study in some detail Case B in \eqref{e:frac}. In
that case, we have
\begin{equation}\label{E:normalizations1}
    G_w (0) = \pm 1,\quad |F_z(0)| = 1,
\end{equation}
for all mappings respecting the normalization
\eqref{E:qamnorm}, and the
possible pairs $(F_z(0),G_w(0))$ are actually restricted
to be in the discrete subset $D(n-\alpha,\mu -1)$ of the torus
${\sf U}(1)\times{\sf U}(1)$
and $D(k,l)$ is defined by
\begin{equation}\label{E:defD}\begin{aligned}
D(k,l)&:=
\left\{
(\gamma,\delta)\in {\sf U}(1)\times\left\{ -1,1 \right\}\colon \gamma^k \delta^l = 1
\right\}\\
&= \left\{ \left(
e^{i \pi \left(\frac{ 2p + ql}{k} \right) }, e^{i \pi q}\right) \colon 1\leq p < k , q = 0,1 \right\}.
\end{aligned}
\end{equation}
  Now assume that
  $M$  satisfies the assumptions of
  Proposition~\ref{P:prelimnorm} with normal coordinates $(z,w)$ chosen such that
  \eqref{E:qamnorm} holds. Then any formal map
  $H\colon (\C^2,0) \to (\C^2,0)$,
  which respects all coefficients $q_{\alpha,\mu}(\chi)$
corresponding to all invariant pairs $(\alpha,\mu)$,
i.e.\ which   satisfies
  $q_{\alpha,\mu} (\chi)= F_z(0)^{\alpha} G_w (0)^{\mu -1}q_{\alpha,\mu} (\bar F(\chi,0))$
  for all $(\alpha,\mu) \in Q_M$, has the property
  \begin{equation}
    (F_z (0), G_w (0)) \quad \in\quad  N_M := \bigcap_{(\alpha,n,\mu) \in \Lambda_M}
    D(n - \alpha, \mu-1).
    \label{e:Nmdef}
  \end{equation}
  In particular, this holds (by Proposition~\ref{P:0basic}) for
  every formal invertible map taking $M$ {\em into itself}.
  (Note, however, that in order for the conclusion \eqref{e:Nmdef} to hold,
  it is not necessary that $H$ maps $M$ into itself, but only that
  the coefficients corresponding to invariant pairs are respected).

  We have the obvious bound for the number of elements
  $$\# N_M \leq  2 \gcd \{ n-\alpha \colon (\alpha,n,\mu) \in \Lambda_M \}=:2 g_M.$$
  This bound is actually sharp, as pointed out in Example~\ref{ex:finauto} above;
  more generally,
  we note that, for any choice of real coefficients $c_{\alpha,n,\mu}$,
  the hypersurface $\tilde M$ given by
\begin{equation}\label{hyp-admissible}
  \imag w = \sum_{ \left( \alpha,n,\mu \right)\in\Lambda_M}
  c_{\alpha,n,\mu} \real(z^\alpha \bar z^{n}) (\real w)^{\mu}
\end{equation}
  has at least as many different automorphisms as $\# N_M$.
In particular, we have the crude bound $\#N_M \leq 2 (n-\alpha)$
as stated in Theorem~\ref{T:finitedet}.

We are now going to study the conditions guaranteeing that the group $N_M$ is trivial.
We first observe that, if all $\mu$'s for $(\alpha,n,\mu)\in \Lambda_M$ are odd,
then any hypersurface in \eqref{hyp-admissible} has the nontrivial automorphism
$(z,w)\mapsto (z,-w)$. Hence we may assume that some of the $\mu$ is even.
Then, to investigate the conditions for the triviality of $N_M$,
let us write $S_k$ for the set of all $k$-th roots of unity; we
  then have
  \[ D(k,l) =
  \begin{cases}
    S_k \times \left\{ 1 \right\} \cup S_k \times \left\{ -1 \right\} & l \text{ even}\\
    S_k \times \left\{ 1 \right\} \cup e^{\frac{i\pi}{k}}S_k\times
    \left\{ -1 \right\}& l \text{ odd}.
  \end{cases}
  \]
  Let us write $e^{\frac{i \pi}{k}} S_k = S_k^-= \{ \gamma \in \mathbb{T} \colon \gamma^k = -1\}$,
  and set
  \[
  g_M^+ := \gcd \left\{ n-\alpha \colon
  (\alpha,n,\mu) \in \Lambda_M , \mu \text{ odd}\right\}, \quad
  g_M^- := \gcd \left\{ n-\alpha \colon
  (\alpha,n,\mu) \in \Lambda_M , \mu \text{ even}\right\},\]
  so that  $\gcd\left\{ g_M^+,g_{M}^{-} \right\} = g_M$.
Here and in the following we use the convention that
$\gcd(\emptyset)=0$ and $S_{0}^- := \emptyset$, $S_{0} := {\sf U}(1)$, and
we continue to use the convention that the $\gcd$ of a family
of numbers possibly containing $0$ is the same as the
$\gcd$ of its nonzero members. With this notation, we have
  \[
  \begin{aligned} N_M &=
    S_{g_M} \times \left\{ 1 \right\} \bigcup
    \left( \bigcap_{\mu \text{ odd} } S_{n-\alpha} \cap \bigcap_{\mu \text{ even}}
    S_{n-\alpha}^-\right)\times \left\{ -1 \right\}\\
    &=  S_{g_M} \times \left\{ 1 \right\} \cup\left(
    S_{g_M^+}\cap \bigcap_{\mu \text{ even}} S_{n-\alpha}^{-}\right)
    \times \left\{ -1 \right\}.
  \end{aligned}.
  \]
  Denoting by $\ord_2 (k)$ the maximal $j$ such that $2^j$ divides $k$, we note that
  \[
  S_k^- \cap S_{\tilde k}^{-} =
  \begin{cases}
    \emptyset & \ord_2 (k) \neq \ord_{2} (\tilde k)
    \\
    S_{\gcd\{k,\tilde k\}}^{-} & \ord_2 (k) = \ord_2 (\tilde k).
  \end{cases}
  \]
Indeed, if $\gamma\in S_k^- \cap S_{\tilde k}^{-}$, i.e.\ $\gamma^k=\gamma^{\tilde k}=-1$, then
$\gamma^{\gcd\{k,\tilde k\}}= \gamma^{ak+\tilde a\tilde k}=\pm1$
and, since $\gcd\{k,\tilde k\}$ divides $k$ and $\tilde k$, we must have
$\gamma\in S_{\gcd\{k,\tilde k\}}^{-}$. On the other hand, if
$\gamma\in S_{\gcd\{k,\tilde k\}}^{-}$, one clearly has
$\gamma^k=\gamma^{\tilde k}=-1$ if and only if $\ord_2 (k) = \ord_2 (\tilde k)$.
  Thus, we have
  \[
  N_M =
  \begin{cases}
    S_{g_M}\times \left\{ 1 \right\} & \ord_2 (n - \alpha) \text{ is not constant for all even } \mu,
    \\
    S_{g_M}\times \left\{ 1 \right\} \cup \left(
    S_{g_M^+} \cap S_{g_M^-}^- \right) \times\left\{ -1 \right\} & \ord_2 (n-\alpha) \text{ is constant
    for all even } \mu
  \end{cases}
  \]

  Similarly $S_{k} \cap S_{\tilde k}^{-}$ (with $k$ possibly zero) is nonempty exactly when
  $\ord_2 k > \ord_2 \tilde k$, in which case $S_k \cap S_{\tilde k}^- =
  S_{\gcd\left\{ k,\tilde k \right\}}^-$, and we obtain
  \begin{equation}
    \# N_M =
    \begin{cases}
      g_M & \ord_2(n-\alpha) \text{ is not constant for all even } \mu  \\
        &\text{or } \ord_2 g_M^- \geq \ord_2 g_M^+ \\
    2 g_M  & \text{otherwise}.
    \end{cases}
    \label{e:nmestimate}
  \end{equation}
  Note that $N_M \cong \Z_{g_M} $ in the first
  case in \eqref{e:nmestimate} and $N_M \cong \Z_{g_M} \oplus \Z_2$ in the second case.
  In particular, if $k$ and $\tilde k$ do not have any common
  divisors, then $S_k\cap S_{\tilde k}^-$ is nonempty (it may only contain
  the point $-1$) if and only if $ k $ is even and $\tilde k$ is odd.
  Thus, if $g_M = 1$ and not all the $\mu$ are odd (as we assumed),
$N_M$ may only contain the two points $(1,1)$ and $(-1,-1)$;
  it contains only the point $(1,1)$
  if and only if in addition to $g_M=1$ we have
$(-1,-1)\notin D(n-\al,\mu-1)$, i.e.\
either $n - \alpha$ is odd for some invariant pair $(\alpha,\mu)$ with $\mu$ odd,
  or $n-\alpha$ is even for some invariant pair $(\alpha,\mu)$
  with $\mu$ even.

\end{rem}

\subsection{The basic equation} In order to construct the normal form, we
need the following technical result.
\begin{prop}\label{P:basicequation}
  Let $M$, $M^\p$ be formal hypersurfaces in $\C^2$, each given in normal
  coordinates at $0\in M$ and $0\in M'$, respectively.  Let
  $H = (F,G)\colon (\bC^2,0)\to (\bC^2,0)$ be a formal invertible map,
  assume that there are
  at least two invariant pairs for $M$ and $M'$, and let $(\al_0,\mu_0)\in Q_M$
  denote one of them.
  Then, for each $k\geq
  0$, there is a  universal polynomial $R_{\al_0,\mu_0}^k$  such that
\begin{multline}\label{E:basicequation}
  \overline{G_w} (0) q_{\al_0,\mu_0+k} (\chi)  =  \\
  \frac{{F_z (0)}^{\al_0} {\overline{G_w} (0)}^{\mu_0}}{k!}  \Biggl[
  \left( \al_0 \frac{F_{zw^k} (0)}{F_z (0)} +
  \frac{\mu_0}{k+1} \frac{\ov{G_{w^{k+1}}}(\chi, 0)}{\ov{G_w (0)}} -
  \frac{\overline{G_{w^{k+1}}} (0)}{\ov{G_w (0)}}  \right) q^\p_{\al_0,\mu_0} (\bar F
  (\chi,0))
  \\
  +  \ov{F_{w^k}} (\chi,0)\,
    q^\p_{\al_0,\mu_0; \chi^\p} (\bar F (\chi,0)) - \ov{F_{w^k}} (0)\,
    q^\p_{\al_0,\mu_0; \chi^\p} (0)  \Biggr] \\
   + R_{\al,\mu}^k \left( q^\p, F_{z^a w^b} (0),F_{w^k} (0),
    \ov{F_{w^l}} (0), \ov{F_{w^l}} (\chi,
    0 ), \ov{G_{w^m}} (0), \ov{G_{w^m}} (\chi,0) \right)
\end{multline}
where $a\leq \al_0$,  $b,l < k$, $m\leq k$, and $q^\p$ is a shorthand notation
for terms of the form ${q^\p}^{(l)}_{\al,\mu} (\bar F(\chi,0))$
and ${q^\p}^{(l)}_{\al,\mu} (0)$.
Moreover,
if $\al_0>1$, then there is, for each $k\geq 0$,  a universal
polynomial $S_{\al_0,\mu_0}^k$ such that
\begin{multline}\label{E:basicequationwk}
  q_{\al_0 - 1, \mu_0 + k}(\chi ) = \frac{1}{k!} q^\p_{\al_0, \mu_0} (\bar
  F(\chi,0)) F_{w^k} (0) F_z(0)^{\al_0 - 1} \overline{G_w} (0)^{\mu_0}\\ +
  S_{\al_0,\mu_0}^k \left( q^\p, F_{z^a w^b} (0), \ov{F_{w^l}} (0), \ov{F_{w^l}} (\chi,
    0 ), \ov{G_{w^m}} (0), \ov{G_{w^m}} (\chi,0) \right)
\end{multline}
where $a\leq \al_0$, $b,l < k$, $m\leq k$.
\end{prop}

\begin{proof} We continue to use the notation from the previous
sections and consider a formal map $H=(F,G)$ sending $M$ into
$M'$. We rewrite
 \eqref{E:start}  as
follows:
\begin{multline}\label{newstart}
  \bar G\big(0,Q(z,\chi,\tau)\big) = \bar G (\chi,\tau) + \sum_{\al,
    \mu} q_{\al, \mu}^\prime \big(\bar F (\chi,\tau)\big)
  \big(F(z,Q(z,\chi,\tau))\big)^\al
  \big(\bar G(\chi,\tau)\big)^\mu \\
  -\sum_{\al, \mu} q_{\al, \mu}^\prime \big(\bar F
  (0,Q(z,\chi,\tau)\big) \big( F(z,Q(z,\chi,\tau))\big)^\al \big( \bar
  G(0,Q(z,\chi,\tau))\big)^\mu.
\end{multline}
As before we will compare the
coefficients of $z^{\al_0} \tau^{\mu_0+k}$ in \eqref{newstart} for
$k\ge 1$. By Taylor expanding and multiplying out, we see that these
coefficients can be written as universal polynomials in
$q_{\al,\mu}(\chi)$, $(q_{\al,\mu}^\prime)^{(l)} (\bar F(\chi, 0))$,
$(q_{\al,\mu}^\prime)^{(l)} (0)$, $\ov{F_{w^l}} (\chi,0)$, $\bar
G_{w^l} (\chi,0)$, $\ov{G_{w^l}} (0)$, and $F_{z^a w^b}(0)$.  We will
put restrictions on $\al$, $\mu$, $l$, $a$, and $b$.

The Taylor expansion of the left hand side yields
\begin{equation}\label{E:Gtaylorexpanded2}
  \bar G (0,Q(z,\chi,\tau)) = \sum_{l=0}^\infty \frac{1}{l!}
  \ov{G_{w^l}} (0) \left(\tau + \sum_{\al, \mu} q_{\al,\mu} (\chi) z^\al
  \tau^\mu\right)^l.
\end{equation}
For $l=1$, the only term with $z^{\al_0} \tau^{\mu_0+k}$ has the
coefficient $\overline{G_w} (0) q_{\al_0,\mu_0+k} (\chi)$.  The general term
coming from
\[
\ov{G_{w^l}} (0) \left(\tau + \sum_{\al, \mu} q_{\al,\mu} (\chi) z^\al
\tau^\mu\right)^l, \quad l \ge 2,
\]
has the form
\begin{equation}\label{E:wastoolong}
  \ov{G_{w^l}} (0) q_{\al_1,\mu_1}(\chi) \cdots q_{\al_m,\mu_m}(\chi) \,
z^{(\sum_{j=1}^{m} \al_j)} \, \tau^{(l-m + \sum_{j=1}^{m} \mu_j)},
  \quad 0\le m\le l,
\end{equation}
where the product and the sum over an empty set of indices are, by
definition, $1$ and $0$ respectively.  For a term with $z^{\al_0}
\tau^{\mu_0+k}$ we must have $\al_j\le\al_0$, which also implies,
since $(\al_0,\mu_0)\in Q_M$,  that $\mu_j \geq \mu_0$. So either
$\al_j< \al_0$ for all $j$ or $\al_1 = \al_0$ and $m = 1$. For
\eqref{E:wastoolong} to be such a term, we must have
\begin{equation}\label{term-rel}
l-m + \sum_{j=1}^{m} \mu_j=\mu_0+k.
\end{equation}
Thus, if $\al_1= \al_0$ we have $l\leq k+1$, and $l = k+1$ yields
the term
\begin{equation}\label{left-term}
\frac{1}{k!} \ov{G_{w^{k+1}}}(0) q_{\al_0,\mu_0}(\chi) z^{\al_0} \tau^{\mu_0+k}.
\end{equation}
For other possible terms, we have $m \ge 1$ and $\mu_j \ge \mu_0+1$ for
all $j$.  In this case \eqref{term-rel} implies $l\le k -
(m-1)\mu_0\le k$ and the resulting coefficient of $z^{\al_0}
\tau^{\mu_0+k}$ will be a polynomial of
\begin{equation}\label{poly-1}
\ov{G_{w^{l}}}(0), \quad l\le k, \text{ and } q_{\al,\mu}, \quad
\al \le \al_0,\quad \mu \le  \mu_0 +k.
\end{equation}
We claim that if  $\mu = \mu_0 + k$, then $\al < \al_0 -1$.
Let us first check that $q_{\al_0,\mu_0+k}$ can not appear.
Indeed, \eqref{term-rel} implies that $l = 1$ in that case, which
we have already separated above. Now, $q_{\al_0 - 1, \mu_0
  +k}$ can only appear if $q_{1,0} \neq 0$, a case excluded by our
  assumption that there are at least two invariant pairs.
Hence, the remaining terms are polynomial
in
\begin{equation}\label{poly-2}
\ov{G_{w^{l}}}(0), \quad l\le k, \text{ and } q_{\al,\mu}, \quad
\al\leq \al_0,\quad \mu < \mu_0 +k \text{ or }  \al < \al_0 - 1, \mu= \mu_0 +k.
\end{equation}

Let us now examine the right hand side of \eqref{newstart}.  Clearly,
the first term on the right does not contribute a term
$z^{\al_0}\tau^{\mu_0+k}$, so let us consider the first sum. If we
get a term containing $\ov{G_{w^l}} (\chi, 0)$, the minimum power
of $\tau$ we get from the other terms is $(\al - \al_0)_+ + \mu -
1$, where we write $n_+ = \max(n,0)$. Hence,
\begin{equation}\label{E:tauexp1}
l + (\al - \al_0)_+ + \mu - 1 \leq \mu_0 +k.
\end{equation}
If $\mu > \mu_0$, \eqref{E:tauexp1} implies that $l<k+1$. On the other
hand, if $\mu\leq \mu_0$, either $\al + \mu > \al_0 + \mu_0$ or
$(\al,\mu) = (\al_0,\mu_0)$. In the first case, we have $\al - \al_0
> \mu_0 - \mu \ge 0$, and so \eqref{E:tauexp1}
again implies that $l<k+1$; in the
second case, $l = k + 1$ can appear, in which case the
corresponding term is given by
\begin{equation}\label{E:rightterm1}
  \frac{\mu_0}{k!} \ov{G_{w^{k+1}}} (\chi, 0) q^\p_{\al_0, \mu_0} (\bar
  F(\chi, 0)) F_z(0)^{\al_0} \overline{G_w} (0)^{\mu_0-1}.
\end{equation}
Any term containing a factor $F_{z^a w^b} (0)$ comes together with
a contribution of $z^a (\tau + \sum q_{\bt, \nu}(\chi) z^\bt
\tau^\nu)^b$. The general term from expanding this product is
\begin{equation}\label{E:genterm2} q_{\bt_1,\nu_1} (\chi) \dots q_{\bt_r, \nu_r} (\chi)
z^{a + \sum
  \bt_j} \tau^{b - r + \sum \nu_j}, \quad 0\leq r\leq b,\end{equation}
where the case $r=0$ is understood to mean $z^a\tau^b$. Hence, we
have $\bt_j < \al_0$ unless for some $j$, $\bt_j = \al_0$, in
which case $r=1$. We will discuss the latter case later and,
hence, we assume for now that $\beta_j<\al_0$ for each
$j=1,\ldots, r$ or that $r=0$. Note that $\bt_j < \al_0$ implies
that $\nu_j
> \mu_0$. Thus, examining the overall exponents of $z$ and $\tau$
in a term from the first sum on the right which contains the
factor $F_{z^aw^b}(0)$ and which contributes
$z^{\al_0}\tau^{\mu_0+k}$, we see that
$$
\mu + b-r+\sum \nu_j + (\al-\al_0-1+a+\sum\beta_j)\leq\mu_0+k
$$
for some $0\leq r\leq b$ and some sequence of indices
$\beta_j<\al_0$, $\nu_j>\mu_0$, $j=1,\ldots, r$ such that
$$
\al-\al_0-1+a+\sum\beta_j\geq 0
$$
and
$$
a+\sum\beta_j-\al_0\leq 0.
$$
Observe that these inequalities imply, under the assumption made
above that $r=0$ or $\beta_j<\al_0$ for $j=1,\ldots r$,
\begin{equation}\label{E:tauwt1}
 \mu + b   + r \mu_0 + \left( \al - \al_0 - 1  + a + \sum \bt_j\right)_+ \leq \mu_0 +
k.\end{equation}
If $\mu > \mu_0$, \eqref{E:tauwt1}  implies $b< k$. On the other hand, if
$\mu \leq \mu_0$, either $(\al, \mu) = (\al_0, \mu_0)$ or $\al + \mu >
\al_0 + \mu_0$. If $(\al, \mu) = (\al_0, \mu_0)$,
  \eqref{E:tauwt1} implies that $b<k$ unless $a+\sum \bt_j \leq 1$.
So either $a = 1$ and $r=0$, which gives us the term
\begin{equation}
  \label{eq:rightterm2}
  \frac{\al_0}{k!}q^\p_{\al_0, \mu_0} (\bar F(\chi, 0)) F_{zw^k}(0)
  F_z(0)^{\al_0 -1} \overline{G_w} (0)^{\mu_0},
\end{equation}
or $a = 0$, $r =1$ and $\bt_1 = 1$; going back to \eqref{E:tauwt1}
we see that $b\le k$, giving rise to terms containing
$F_{w^{k}} (0)$.
 If $\al + \mu > \al_0 + \mu_0$, \eqref{E:tauwt1} implies $a+b
\leq k$.

We are now going to check which $q_{\bt,\nu}$ can appear from the
first sum.
If in \eqref{E:genterm2} one of the $\bt_j$, say $\bt_1$ satisfies
$\bt_1 = \al_0 - 1$ and the corresponding $\nu_1 = \mu_0 +k$, we
clearly have $r\leq 2$. If $r = 2$, then $\beta_2=1$, $a=0$,
$\nu_2 = 0$, $\mu = 0$. Since $(\beta_2,\nu_2)=(1,0)$, terms of
this form cannot appear (by our assumption that there are at least two
invariant pairs).
 On the other hand, if $r = 1$, we have $a \leq 1$,
 and also $\al + a \leq 2$. Checking the possible cases $\al=1,2$
 and $a=0,1$, we see that all these are again excluded by the assumption
 of having at least two  invariant pairs.

Let us now turn to the case where, for some $j$, $\al_0 = \bt_j$.
It follows immediately that $r=1$, $a=0$, and $\nu_1 \geq \mu_0$,
and instead of \eqref{E:tauwt1} we now obtain
\begin{equation}
  \label{eq:tauwt1a}
  \mu + b-1 + \nu_1 + (\al  - 1)_+ \leq \mu_0 + k.
\end{equation}
If $\nu_1 > \mu_0$, we immediately obtain $b \leq k$. If on the other
hand, $\nu_1 = \mu_0$, we have $b \leq k$ unless $\mu = 0$. If $\mu =
0$, again $b\leq k$, this time unless $\al= 1$. It follows that a term
with $b = k+1$ can only appear if $(1,0)$ is an invariant pair, which
is impossible by assumption.

We also note that $\nu_1 = \mu_0 +k$ implies $\mu + b + (\al - 1)_+
\leq 1$. Since in that case, $b \geq 1$, this implies that $(1,0)$ is
an invariant pair, which
is excluded.

At this point, let us note that for the $q_{\beta,\nu} (\chi)$
coming from the first sum, we have shown that either $\nu<\mu_0 +
k$ or $\beta < \al_0 - 1$ if $\nu = \mu_0 + k$.

We now turn to the expansion of the term $q^\p_{\al,\mu} (\bar
F(\chi,\tau))$, which appears as a factor in the first sum on the
right in \eqref{newstart}, as a Taylor series in $\tau$. The
coefficient of $\tau^l$ is a linear combination of terms of the
form
\[ {q^\p_{\al,\mu}}^{(p)}(\bar F(\chi,0))  \ov{F_{w}} (\chi,0)^{p_1}
\dots \ov{F_{w^k}} (\chi,0)^{p_k }, \] where $p_1 + 2 p_2 + \dots
k p_k = l$. If this term appears as a factor in a coefficient of
$z^{\al_0}\tau^{\mu_0+k}$, it follows that $l + \mu + \max(\al -
\al_0 , 0) \leq \mu_0 + k$, and by arguments with which the patient
reader
is familiar by now, either $l< k$, or $\al = \al_0$, $\mu =
\mu_0$, and $l=k$. We conclude that the only term containing $\bar
F_{w^k}(\chi,0)$ is
\begin{equation}
  \label{eq:lastterm}
  q^\p_{\al_0,\mu_0; \chi^\p} (\bar F (\chi,0)) \ov{F_{w^{k}}} (\chi,0).
\end{equation}

We are left with contributions from the second sum on the right hand
side of \eqref{newstart}. We have already expanded $\bar G
(0,Q(z,\chi, \tau))$ in \eqref{E:Gtaylorexpanded}. The general term
here has the form \eqref{E:wastoolong}. Hence, $\sum \al_j \leq
\al_0$, which implies that either $\al_j = \al_0$ for some $j$ (in
which case, $m = 1$), or that $\al_j < \al_0$ for all $j$. If $\al_j =
\al_0$, there must be at least one $\tau$ contributing from the first
factor of the summand
\[q_{\al, \mu}^\prime \big(\bar F (0,Q(z,\chi,\tau)\big) \big( F(z,Q(z,\chi,\tau))\big)^\al
\big( \bar G(0,Q(z,\chi,\tau))\big)^\mu;\] we then see that $\mu +
l - 1 + \mu_1 + \al_0 \leq \mu_0 +k$. But $\al_0 \geq 1$, so $\mu
+ \mu_1 + l \leq \mu_0 + k$; since $\al_1 = \al_0$, it follows
that $\mu_1 \geq \mu_0$, and so $l\leq k$. We also note that since
$l\geq 1$, then $\mu_1 < \mu_0 +k$.  On the other hand, if $\al_j
< \al_0$, $\mu_j
> \mu_0$ for all $j$, the first term contributes either a term
containing a $\tau$, or (only if $\mu_0 = 0$) it contributes a
term containing $z^{\al_0}$; in the first case, we have $\mu + l +
m \mu_0 + \max(0, \al - \al_0) \leq \mu_0 +k$, from which it
follows that $l\leq k$. In the second case, we must have $m = 0$
and $\mu - 1 + l + \al \leq \mu_0 +k$, from which we again
conclude that $l\leq k$ (since $\al_0 \geq 1$).

A term containing an $F_{z^a w^b}(0)$ comes with a term of the form
\eqref{E:genterm2}. Again, the first factor $q_{\al,\mu}^\p (\bar
F(0,Q(z,\chi,\tau))$ contributes either a $\tau$ or a term
$z^{\al_0}$. In the first case, if it contributes at least a $\tau$,
instead of \eqref{E:tauwt1}, we get
\begin{equation}
  \label{eq:tauwt2}
  \mu + b + 1 + r\mu_0 + \max\left(0,\al - \al_0 - 1 + a + \sum
  \bt_j\right) \leq
  \mu_0 +k,
\end{equation}
if $\bt_j < \al_0$ for all $j$. If $\mu \geq \mu_0$, \eqref{eq:tauwt2}
implies $b<k$. On the other hand, if $\mu<\mu_0$, $\al_0 + \mu_0 > \al
+ \mu$, and so \eqref{eq:tauwt2} again implies $b< k$. If $\bt_j =
\al_0$ for some $j$, $r=1$ and $a=0$, and instead of \eqref{eq:tauwt2}
we obtain
\begin{equation}
\label{eq:tauwt2a} \mu + b + \nu_1 +  \al - 1 \leq \mu_0 + k.
\end{equation}
Since $\al_0 = \bt_1$, we have $\nu_1 \geq \mu_0$. Hence,
\eqref{eq:tauwt2a} implies $b \leq k$, and if $b = k$ in
\eqref{eq:tauwt2a} we have $\mu + \al - 1 \leq \mu_0 -\nu_1 \leq 0$.
Hence, this implies that $(1,0)$ is the invariant pair, which
is excluded.
In the second case, if we get a $z^{\al_0}$ from $q_{\al,\mu}^\p (\bar
F(0,Q(z,\chi,\tau))$, we must have $a = r = 0$, and
\begin{equation}
  \label{eq:tauwt3}
  \mu + b + \al \leq \mu_0 +k,
\end{equation}
which again implies $b <k$. We also note that any $q_{\al,\mu}$ coming
from here satisfies $\mu < \mu_0 + k$.

We now come to contributions from expanding $q_{\al,\mu}^\p (\bar
F(0,Q(z,\chi,\tau))$. A term containing an $\ov{F_{w^l}} (0)$ comes
with a term of the form \eqref{E:wastoolong} with $G$ replaced by $F$.
Hence, we again either have $\al_j <\al_0$ for all $j$, or $\al_1 =
\al_0$ and $m =1$. In the first case, $\mu_j > \mu_0$ for all $j$, and
so
\begin{equation}\label{eq:tauwt4}
 \mu + l + m \mu_0 + \max\left({0 , \al - \al_0 + \sum \al_j}\right)  \leq \mu_0
+k. \end{equation}
If $\mu > \mu_0$, this implies $l<k$. If on the other hand, $\mu \leq
\mu_0$, we have either $\al + \mu > \al_0 + \mu_0$ or $(\al,\mu)=
(\al_0, \mu_0)$. If $\al + \mu > \al_0 + \mu_0$, $\al > \al_0$, and so
\eqref{eq:tauwt4} becomes
\[\mu + l + m \mu_0 + \al - \al_0 + \sum \al_j   \leq \mu_0
+k, \] from which it follows that $l<k$.  In fact, if $(\al,\mu)=
(\al_0, \mu_0)$, \eqref{eq:tauwt4} turns into $l + m\mu_0+\sum \al_j
\leq k$, so again $l\le k$.

In the second case, if $\al_1 = \al_0$, we have $\mu_1 \geq \mu_0$,
and $\mu_1 + l - 1 + \al + \mu \leq \mu_0 +k$. Unless $\mu =0$ and
$\al =1$, this implies that $l<k$.

Similarly as before, the $q_{\al,\mu}$ from this first factor satisfy either
$\mu < \mu_0 + k$ or $\al < \al_0 -1$.

By now, we have nearly finished the proof of \eqref{E:basicequation}.
The only difference is that we have constructed a remainder term $\td
R$ which depends also on $q_{\bt,\mu}$, where either $\bt < \al_0 - 1$
or $\bt = \al_0$, but $\mu<\mu_0 +k$. The proof will be finished by
showing that terms of the form $q_{\bt,\mu}$ with $\bt < \al_0$
themselves depend on terms which appear in the remainder and on
$q_{\gm,\nu}$ where $\nu < \mu_0 +k$. Applying this result repeatedly
then gives us the form of the remainder we seek. We incorporate the
proof of this auxiliary statement in the derivation of
\eqref{E:basicequationwk}, which follows.

We shall extract the coefficient of $z^{\bt} \tau^{\mu_0 + k}$,
where $\bt<\al_0$, on both sides of \eqref{newstart}. Let us start
with the left hand side. First, we get the term on the left hand
side of \eqref{E:basicequationwk} from the Taylor expansion
\eqref{E:Gtaylorexpanded2} with $l = 1$. A term of the form
\eqref{E:wastoolong} can only appear if $\al_j \leq \bt < \al$;
furthermore, we have $m > 0$. Now $\al_j < \al$ implies that
$\mu_j > \mu \geq 0$ for all $j$, so that $\mu + k = \sum_{j=1}^m
\mu_j + l - m \geq \mu + l$, or $l\leq k$. Furthermore, all terms
$q_{\gm,\nu}$ coming from \eqref{E:Gtaylorexpanded2} other than
$q_{\bt, \mu_0 +k}$ satisfy $\nu < \mu_0 +k$, as can be easily
checked using \eqref{term-rel}.

Now, we turn to the right hand side of \eqref{newstart}. Of course,
the first term does not contribute as long as $\al_0 > 1$; for the
first sum, as usual, we distinguish the cases $\mu \leq \mu_0$ and
$\mu > \mu_0$.

We first discuss the case $\mu \leq \mu_0$. If a derivative $\bar
G_{\tau^l}(\chi,0)$ contributes, then from the other terms in the
product, the minimum exponent of the $\tau$'s is $\mu - 1 + (\al -
\bt)_+$. Since $(\al_0, \mu_0)$ is an invariant pair, $\al - \bt > \al
- \al_0 \geq \mu_0 - \mu \geq 0$. So the overall exponent of $\tau$ is
at least $\mu + l - 1 + \al - \bt$. It follows that
\[ \mu + l + \al - \al_0 \leq \mu_0 +k, \]
from which we conclude $l\leq k$.  If we have a term containing
$F_{z^a w^b}(0)$, using the notation from \eqref{E:genterm2} the other
terms contribute at least
\[ \mu + \left(\al - 1 - \bt + a + \sum \bt_j\right)_+\]  to the overall exponent of
$\tau$; as before, $\al - \bt > \al - \al_0$, and since $\bt_j \leq
\bt < \al_0$, $\mu_j > \mu_0$, so we have
\begin{equation}\label{eq:tauwt5} \mu  + b + r\mu_0 + \al - \al_0 + a + \sum\bt_j \leq
\mu + b +r \mu_0 + \al - \bt - 1 + a + \sum \bt_j \leq
\mu_0 + k\end{equation} and therefore $b\leq k$. If $b = k$, letting $D = \mu - \mu_0 + \al - \bt -1 \geq 0$, we get
\[  D + r\mu_0 + a + \sum\bt_j \leq 0, \]
so we must have $D=r \mu_0 = a = \bt_j = 0$ for all $j$, which in
particular implies $r = 0$. Furthermore, since $D = 0$, and $\mu -
\mu_0 + \al-\bt > \mu - \mu_0 +\al - \al_0 > 0$ if $(\al,\mu) \neq
(\al_0,\mu_0)$, it follows that $(\al,\mu) = (\al_0,\mu_0)$ and $\bt =
\al_0 -1$, which leads to the term
\begin{equation}
\label{eq:bequfwk1} \frac{1}{k!} q^\p_{\al_0, \mu_0} (\bar F(\chi,0)) F_{w^k} (0)
F_z(0)^{\al_0 - 1} \overline{G_w} (0)^{\mu_0}
\end{equation}

If a term $\ov{F_{w^l}} (\chi, 0)$ from expanding the first factor
appears, the exponent of $\tau$ is at least $l + \al - \bt + \mu$; a
similar discussion as above shows that this implies $l<k$.

We now turn to the case $\mu > \mu_0$. In this case, a derivative
$\bar G_{\tau^l} (\chi, 0)$ will contribute only if $\al - \bt + l
+\mu - 1 \leq \mu+ k$ which implies $l\leq k$.  Now let us consider
the derivatives $F_{z^a w^b} (0)$, for which we have
\eqref{eq:tauwt5}, which in our case now implies $q < k$.  Also, as
above, a term $\ov{F_{w^l}} (\chi,0) $ only appears if $l<k$.

Now we have to check the second sum on the right hand side. A
discussion similar to the one above, using the fact that any term
from the first factor in each of the products contributes at least
one power of $\tau$, now shows that only terms of the form claimed
appear, which finishes the proof of \eqref{E:basicequationwk}, and
as explained before, also that of \eqref{E:basicequation}.
\end{proof}

\subsection{Construction of the normalization map} Recall that we
assume that $\Lambda_{M^\prime}$ consists of at least two points
$(\al,\mu)\neq(\td \al,\td \mu)$. We also assume that one of these
pairs, say $(\al,\mu)$, satisfies $\al \neq n$, where
$n=n(\al,\mu)$ as defined by \eqref{E:nam}.  We also note that at
least one of $\al$ and $\al^\p$ is greater than $1$; let us denote
the corresponding invariant pair by $( \al^\p, \mu^\p)$; thus,
$\al'>1$.

We first normalize $M$ as in Proposition \ref{P:prelimnorm}. Hence
we shall assume that \eqref{E:qamnorm} holds and
\begin{equation}\label{normnow}
F(z,0) = F_z (0) z, \quad
(F_z(0), G_w (0))\in C(\alpha,n,\mu) \cap C(\tilde{\alpha}, \tilde{n}, \tilde{\mu}),
\end{equation}
where $C$ is given by \eqref{e:defC}.
We shall construct a map $H= (F,G)$ with these first order
derivatives $F_z(0)$ and $G_w(0)$ given such that the equation of $M = H^{-1} (M^\p)$ has a
special form. The map $H$ will be unique up to at most one more
parameter. The precise formulation of the result will be given in
Theorem \ref{T:normalform} below.

Observe that $G_w(z,0)$ and $F(z,0)$  are uniquely determined by
$F_z(0)$ and $G_w(0)$ in view of Lemma \ref{L:zwnormalifandonlyif}
(note that there $m_0\geq 2$, since $(1,0)$ is not an invariant pair)
and Proposition \ref{P:prelimnorm}. Let $k\geq 1$ and assume that
we have constructed $\ov{F_{w^l}} (\chi,0)$ and $\bar
G_{w^{l+1}}(\chi,0)$ for $l<k$. We first use
\eqref{E:basicequationwk} for $(\al^\p, \mu^\p)$ to determine
$F_{w^k}(0)$ uniquely by the new requirement that $q_{\al^\p-1,
\mu^\p +k} (\chi)$ does not contain any term $\chi^{ n^\prime}$. Let us
now rewrite the basic equation \eqref{E:basicequation} with the
normalizations we have made so far:
\begin{multline}\label{E:basicequation2}
  \overline{G_w} (0) q_{\al,\mu+k} (\chi) = \frac{2i \overline{G_w}(0)}{ k!}
  \left(\frac{n}{ \overline{F_z} (0)} \ov{F_{w^k}} (\chi,0) \chi^{n-1}
    +
    \frac{\al}{F_z (0)} F_{zw^k} (0) \chi^{n}\right) \\
  + \frac{2i}{k!}\Bigl(\frac{\mu}{k+1} \ov{G_{w^{k+1}}} (\chi,0)
  -\ov{G_{w^{k+1}}} (0)\Bigr) \chi^{n}\\
  + \td R_{\al,\mu}^k \left(q^\p, F_{z^a w^b} (0), F_{w^k} (0),\ov{F_{w^{l}}} (0),
    \ov{F_{w^{l}}} (\chi,0) ,\ov{G_{w^{m}}} (0), \bar
    G_{w^{m}} (\chi,0)\right),
\end{multline}
where $\td R_{\al,\mu}^k$ contains the same terms as $R_{\al,\mu}^k$ in \eqref{E:basicequation},
where $F_{w^k} (0)$ has already been determined.
We note that by Lemma~\ref{L:zwnormalifandonlyif}, if $(z,w)$ are
normal coordinates for $M$, then $\ov{G_{w^{k+1}}} (\chi,0)$ is
determined by $\real \ov{G_{w^{k+1}}}(0) =: s_{k+1}$, the
previously determined derivatives, and $Q^\prime$. For
$s_{k+1}\in\R$ we define $\ov{G_{w^{k+1}}}(\chi,0)$ by
\eqref{E:zwnormalifandonlyif} and rewrite \eqref{E:basicequation2}
as
\begin{multline}\label{E:basicequation3}
  \overline{G_w} (0) q_{\al,\mu+k} (\chi) = \frac{2i \overline{G_w}(0)}{ k!}
  \left( \frac{n}{\overline{F_z} (0)} \ov{F_{w^k}} (\chi,0) \chi^{n-1}
    + \frac{\al}{F_z (0)} F_{zw^k} (0)\chi^{n}\right)
  \\+  \frac{2i}{k!}\Bigl(\frac{\mu}{k+1} -1 \Bigr) s_{k+1} \chi^{n}\\
  + \td R_{\al,\mu}^k \left(Q^\p, F_{z^a w^b} (0), F_{w^k} (0), \ov{F_{w^l}}
    (0), \ov{F_{w^{l}}} (\chi,0) ,\ov{G_{w^{m}}} (0), \bar
    G_{w^{m}} (\chi,0)\right),
\end{multline}
where we still write $\td R_{\al,\mu}^k$ for the remainder, which is again a
universal polynomial. Now we observe that we can uniquely
determine $\ov{F_{z^a w^k}} (0)$ for $a>1$ by requiring that
$q_{\al,\mu+k}(\chi)$ is a polynomial of degree at most $n$. We
are now going to examine the coefficient of $\chi^n$ on the right
hand side of \eqref{E:basicequation3}. It is
\begin{equation}\label{E:basic-coefficient}
\frac{2i}{k!}\left( n \frac{\ov{F_{z w^k}} (0)\overline{G_w} (0)}{\bar
F_{z}(0) } + \al \frac{F_{zw^k} (0)\overline{G_w} (0)}{F_z (0) }
+\left(\frac{\mu}{k+1} - 1\right) s_{k+1} \right) +
A_{\al,\mu}^{k,n},
\end{equation}
where $A_{\al,\mu}^{k,n}$ denotes the coefficient of $\chi^n$ in
$\td R_{\al,\mu}^k(\dots )$. Writing
\begin{equation}\label{E:replaceFzwk}
  \frac{F_{zw^k} (0)\overline{G_w} (0)}{F_z (0) }  = B_k
\end{equation}
and noting that $\overline{G_w}(0)$ is real,  we observe that in order
to make the coefficient of $\chi^n$ in $q_{\al,\mu+k} (\chi)$
vanish, we need to solve the equation
\begin{equation}\label{E:equationforBk1}
  n \bar B_k + \al B_k + \left(\frac{\mu}{k+1} - 1\right) s_{k+1} = A_k
\end{equation}
for some right hand side $A_k$. Since $n\neq \al$, this can always
be done, although $B_k$ and $s_{k+1}$ are not uniquely determined.
In order to determine $B_k$ (and hence $F_{zw^k}(0)$ in terms of the first jet of $H$)
uniquely we need to use our second invariant pair $(\td \al,
\td \mu)$. In the basic equation for $q_{\td \al,\td \mu + k} (\chi)$
we see that the coefficient of $\chi^{\td n}$ on the right hand
side has a similar form to the corresponding one of $\chi^n$ in
\eqref{E:basicequation3} above. In fact, inspecting
\eqref{E:basicequation} we see that this coefficient has the form
\begin{equation}\label{E:basic-coefficient2}
  \frac{ ic\varepsilon }{k!}\left( \td n
\frac{\ov{F_{zw^k}} (0)\overline{G_w} (0)}{\ov{F_{
z}}(0) } + \td\al \frac{F_{zw^k} (0)\ov{ G_w} (0)}{F_z (0) }
+\left(\frac{\td \mu}{k+1} - 1\right) s_{k+1} \right) +
A_{\al^\p,\mu^\p}^{k,n},
\end{equation}
where $c>0$ and $\varepsilon$ are the invariant(s) given by
Proposition \ref{P:prelimnorm}.  So if we can solve the real part
of the equation
\begin{equation}\label{E:equationforBk2}
  \td n \bar B_k + \td \al B_k + \left(\frac{\td \mu}{k+1} - 1\right) s_{k+1}
  = A_k^\p,
\end{equation}
then we can make the real part of $i\varepsilon^{-1}$ times the
coefficient of $\chi^{\td n}$ in $q_{\td \al, \td \mu} (\chi)$
vanish, i.e.\ the coefficient
of $\chi^{\td n}$ in the expansion of $q_{\td \al,\td \mu+k}(\chi)$
is $r\varepsilon$ with $r\in\bR$.
This is the additional normalization condition that allows us to solve uniquely for $B_k$ and
$s_{k+1}$, and hence for $F_{zw^k} (0)$ and $\real \bar
G_{w^{k+1}} (0)$. Indeed, if we write $B_k = a_k + i b_k$ and
separate \eqref{E:equationforBk1} and \eqref{E:equationforBk2}
into real and imaginary parts,  then we obtain the following
system of real equations:
\begin{equation}\label{E:realequationsforBkandsk+1}
\begin{array}{llll}
(\al + n) a_k &
    &+ \left(\frac{\mu}{k+1} - 1\right) s_{k+1}
    &=r_k^1\\
&(\al - n) b_k
    & & = r_k^2 \\
(\td \al + \td n) a_k &
    &+ \left(\frac{\td \mu}{k+1} - 1\right) s_{k+1}
    &=r_k^3.
\end{array}
\end{equation}
Since $\al - n$ is not zero, we can solve for $b_k$. Let us
consider the  determinant $D_k$ of the remaining equations in
\eqref{E:realequationsforBkandsk+1} for $a_k$ and $s_{k+1}$,
\begin{equation}\label{E:Jacforaands}
\begin{aligned}
  D_k = &\, \det
  \begin{pmatrix} \al + n &  \frac{\mu}{k+1} - 1\\
    \td \al + \td n & \frac{\td \mu}{k+1} - 1
    \end{pmatrix}\\
    =&\,
    \frac{((\td \alpha+\td n)-(\alpha+n))(k+1)+(\alpha+n)\td\mu-
    (\td\alpha+\td n)\mu}
    {k+1}.
    \end{aligned}
\end{equation} Observe (since $\td\mu\neq\mu$) that $D_k\neq 0$ for $k\geq1$ unless
$\alpha+n\neq \td\al+\td n$ and the number
\begin{equation}\label{E:lambda}
\lambda=\lambda(\alpha,\mu,n,\td\al,\td\mu,\td n):=\frac{(\td\al+\td n)\mu-
(\alpha+n)\td\mu}{(\td\al+\td n)-(\alpha+n)}-1
\end{equation}
is an integer $\geq 1$. In the latter case, we shall allow $s_{\lambda+1}$
to be a real parameter and solve the second equation in
\eqref{E:realequationsforBkandsk+1} for $a_\lambda$ in terms of
$s_{\lambda+1}$. If $\alpha+n=\td\al+\td n$ (which, in particular, implies
that one of the pairs, say $(\alpha,\mu)$, satisfies the assumption
$\alpha\neq n$ made at the beginning of this section) or $\lambda$ is
not an integer $\geq1$, then we may solve the first three equations in
\eqref{E:realequationsforBkandsk+1} uniquely for $a_k$, $b_k$, and
$s_{k+1}$ for every $k$; the normalization map is then determined by
$F_z (0)$ and $G_w(0)$.

Let us now summarize the normalization conditions. Let $M$ be a
formal manifold, given in some fixed system of normal coordinates
$(z,w)$ by $w = Q(z,\chi,\tau)$. We choose an invariant pair
$(\al,\mu)$ with $\al\neq n(\al,\mu) = n$, and an additional
invariant pair $(\td\al,\td \mu)$. We let $(\al^\p, \mu^\p)$
denote $(\al, \mu)$ or $(\td\al ,\td\mu)$ in such a way that
$ \al^\p>1$. Let $K$ denote the set of integers $k\geq 1$ such
that $D_k$ in \eqref{E:Jacforaands} is zero; thus, $K$ is either
empty (if $\alpha+n=\alpha'+n'$ or $\lambda$, given by
\eqref{E:lambda}, is not an integer $\geq 1$) or consists of one
point, namely $\lambda$. We shall say that $Q$ is in its normal form
if the following hold:
\begin{equation}\label{E:normalform1}
    q_{\al,\mu} (\chi ) = 2i \chi^n,\quad q_{\td\al,\td\mu} (\chi ) =
    2ic\varepsilon\chi^{\td n} + O(\chi^{\td n + 1}),
\end{equation}
\begin{equation}\label{E:normalform2}
  q_{\al,\mu+k} (\chi) \text{ is a polynomial of degree at most }
  n-1
  \text{ for } k\geq 1,
\end{equation}
\begin{equation}
  \label{E:normalform2a}
  q_{\al'- 1, \mu' +k}^{(n')} (0) = 0, \text{ for } k \geq 1,
\end{equation}
and
\begin{equation}\label{E:normalform4}
  \imag \varepsilon^{-1} q_{\td\al,\td \mu+k}^{(\td n)} (0) = 0\
  \text{{\rm for $k\geq 1$ and
  $k\not\in K$}},
\end{equation}
where $c$ and $\varepsilon$ are the invariants given by
Proposition \ref{P:prelimnorm}. We have proved the following:

\begin{thm}\label{T:normalform}
  Let $M\subset\C^2$ be a formal hypersurface such that $Q_M$ contains
  at least two points $(\alpha,\mu)\ne (\tilde{\alpha},\tilde{\mu})$
  and (say) $(\al,\mu)$
  satisfies $\al \neq n(\al,\mu)$. Then there exists a system of
  normal coordinates $(z,w)$ for $M$ such that $M$ is given by $w =
  Q(z,\chi,\tau)$ and $Q$ is in its normal form, that is, $Q$ satisfies
  \eqref{E:normalform1}-\eqref{E:normalform4}.  Furthermore, if
  $(z,w)$ are such coordinates, then for any other such system of
  coordinates $(z^\p, w^\p)$, where $z = F(z^\p,w^\p)$, $w = G(z^\p,
  w^\p)$, $F$ and $G$ are uniquely determined by $F_z (0)$, $G_w (0)$
  and $\real G_{w^{k+1}} (0)$ for $k\in K$ and $(F_z (0),G_w (0))$
  belongs to the set $C(\alpha,n,\mu)\cap
  C(\tilde{\alpha}, \tilde{n} ,\tilde{\mu} )$
  (where $C$ is defined in  \eqref{e:defC}). Here, $K$ is the set
  (consisting of at most one point) of integers
  $k\geq 1$ such that $D_k=0$ in
  \eqref{E:Jacforaands}.
  The set $K$ is empty if and only if $n +\al = \td n +
  \td \al$ or $\lambda$ given by \eqref{E:lambda} is not an integer
  $\geq 1$.
\end{thm}

We mention that we have also proved the following Jet
Parametrization Theorem. We write $G_0^k (\C^2)$ for the jet group
of order $k$ of biholomorphic mappings of neighborhoods of $0$ in
$\C^2$ fixing $0$. This theorem implies that $\Aut (M,0)$ can be
embedded as closed real subgroup of $G_0^k (\C^2)$
and hence as its Lie subgroup, for some $k$.

\begin{thm}\label{T:jetparam}
  Let $M\subset\C^2$ be a formal hypersurface such that $Q_M$ contains
  at least two points $(\alpha,\mu)\ne (\tilde{\alpha},\tilde{\mu})$
  and (say) $(\al,\mu)$
  satisfies $\al \neq n(\al,\mu)$. Then there exists an integer $k$ and a formal
  power series map $\psi(Z,W) = \sum_{\alpha \in \N^2} \psi_\alpha (W) Z^{\alpha}$,
  $Z\in\C^2$, with rational coefficients $\psi_\alpha (W)$ in
$W\in G_0^k (\C^2)$ with no poles in $G_0^k (\C^2)$, such that
  \begin{equation} \label{e:param}
    H(Z) = \psi (Z, j_0^k H)
  \end{equation}
    for any $H\in \Aut(M,0)$.
  Furthermore $k$ can be chosen to be $1$ if $n + \alpha = \tilde{n} +
  \tilde{\alpha}$
  or
  $ \frac{ \mu (\tilde{n} + \tilde{\alpha}) - \tilde{\mu} (n+\alpha)}{
  (\tilde{n} + \tilde{\alpha}) - (n + \alpha)} \notin \N_{\ge2}$; otherwise,
  $k$ can be chosen to be  $ \frac{ \mu (\tilde{n} + \tilde{\alpha}) - \tilde{\mu} (n+\alpha)}{
  (\tilde{n} + \tilde{\alpha}) - (n + \alpha)}\in \N_{\ge2}$.
\end{thm}

It is worthwhile to note that in the normal form described above, if
$K$ is the empty set (this in particular means we are in Case B in
\eqref{e:frac}), the normalization group $C(\alpha,n,\mu) \cap
C(\tilde{\alpha}, \tilde{n},\tilde{\mu})$ is {\em discrete}, and it
acts on the space of normal forms by the linear transformations
$(z,w) \mapsto (\gamma z,\delta w)$, $(\gamma,\delta)\in
C(\alpha,n,\mu) \cap C(\tilde{\alpha} , \tilde{n},\tilde{\mu})=
D(n-\alpha,\mu-1) \cap D(\tilde{n} -\tilde{\alpha},\tilde{\mu} -
1)$, where $D$ is defined by \eqref{E:defD}.
Thus, the normal form described above gives (in this case) a
complete solution to the equivalence problem, and it also linearizes
the action of the automorphism group $\Aut (M,p)$. In particular, if
$D(n-\alpha,\mu -1) \cap D(\tilde{n} - \tilde{\alpha}, \tilde{\mu}
-1) = \left\{ (1,1) \right\}$, the power series coefficients of $Q$
(in the normal form) form a complete set of biholomorphic invariants of $M$.
Following the arguments of Remark~\ref{R:prelimnorm} we obtain:

\begin{cor} Let $\Lambda \subset \N^3$ contain two points
  $(\alpha,n,\mu)\ne (\tilde{\alpha},\tilde{n},\tilde{\mu})$ with $n\ne \alpha$, such
  that either
  $\alpha + n= \tilde \alpha+\tilde n $ or
  $\frac{\tilde{\mu} (\alpha+ n) -
  \mu(\tilde{\alpha} + \tilde{n})}{(\alpha+ n) -
  (\tilde{\alpha} + \tilde{n})}$ is not a positive integer. Assume in
  addition that furthermore
  $\gcd\left\{ n-\alpha,\tilde n - \tilde \alpha \right\} = 1$,
  either $\mu$ or $\tilde{\mu}$ is odd, and
  either $\mu + n - \alpha$ or $\tilde \mu + \tilde{n} - \tilde{\alpha}$
  is even. Then two hypersurfaces $(M,p)$ and $(M',p)$ satisfying
  $\Lambda_M = \Lambda_{M'} = \Lambda$ are biholomorphically
  equivalent if and only if their normal forms
  $w = Q(z,\bar{z},\bar{w})$ and $w = Q^{\prime}(z,\bar{z},\bar{w})$ coincide,
i.e.\   $Q\equiv Q'$.
\end{cor}

\section{Dependence on higher order Jets: an Example}
\label{S:Example}

It is natural to ask whether or not the maps {\em really} depend on
higher order jets if the set $K$ is nonempty. In this section we
will give an example of a hypersurface $M$ which satisfies the
assumptions of Theorem~\ref{T:normalform}, the set $K$ consists
exactly of one integer $\ell+1$, and the biholomorphisms are
determined by their $(\ell+1)$-jet, but not their $\ell$-jets at $0$,
where $\ell\ge 5$.
For this, choose positive integers $a,b,c,d$
satisfying
\[ c ( \ell - b) = a (\ell - d), \quad c < a + b - d, \quad \ell/2 < b<d<\ell, \quad a,c>0,
\]
implying $c<a$. For instance, we can take $a=4$, $c=2$, $b=\ell-2$, $d=\ell-1$
 to satisfy all the inequalities.

We then consider the preimage of the quadric $S$, given by  $\imag \eta =
|\zeta|^2$, under the map $B\colon(z,w)\mapsto (z^a w^b + z^c w^d,
w^\ell)$. We first claim that $B^{-1} (S)$ contains a unique
real hypersurface $M$ of the form
\[ t = s^{2b - \ell +1} \phi( z, \bar z, s), \quad w = s+it, \]
where
\begin{equation}\label{desired}
\phi(z,0,s)=\phi(0,\chi,s)=0,\quad \phi (z,\bar z,s)
=z^a\bar z^a+ 2 \real z^c\bar z^a s^{d-b} + O(s^{d-b+1}),
\end{equation}
 and that for this hypersurface $M$,
$\Lambda_M$ contains the points $(a,a,2b - \ell +1)$ and
$(c,a,b + d - \ell + 1)$. Rewriting the equation for $B^{-1} (S)$, we see that
this set is given by
\begin{equation}\label{unnamed}
\sum_{0\le j\le [\frac{\ell-1}2 ]} \binom{\ell}{2j+1} (-1)^{j} s^{\ell - 2j -1} t^{2j+1} = (s^2 + t^2)^b
\left( |z|^{2a} + 2 \real z^c \bar z^a w^{d-b} + |z|^{2c} |w|^{2(d-b)} \right).
\end{equation}
Substituting $t = s^m \lambda$, we obtain
\begin{multline} \label{unnamed1}
\sum_{0\le j\le [\frac{\ell-1}2 ]} \binom{\ell}{2j+1} (-1)^{j} s^{\ell +(m -1)(2j +1)} \lambda^{2j+1} \\ = s^{2b}(1 + s^{2(m - 1)}\lambda^2 )^b
\left( |z|^{2a} + 2 \real z^c \bar z^a (s + i s^m \lambda)^{d-b} + |z|^{2c} |s+i s^m \lambda|^{2(d-b)} \right).
\end{multline}
Note that by this substitution we may only loose the set of solutions of \eqref{unnamed}
corresponding to $s=0$, which does not contain any real hypersurface.

If we set $m:=2b-\ell +1>1$, then we can divide \eqref{unnamed1} by
$s^{2b}$ and get \begin{equation}\label{E:phi}
 \ell \lambda   =
|z|^{2a} + F(z,\bar z, \lambda,s), \end{equation}
 where $F(z,\bar z,\lambda,s)$ is a real-analytic real-valued function that satisfies
\begin{equation}\label{E:F} F(z,\chi,\lambda,0) =
F(z,0,0,s)=F(0,\chi,0,s)=0\quad F_\lambda(0,0,0,0) = 0.
\end{equation} Thus, we may use the Implicit Function Theorem to
solve \eqref{E:phi} for $\lambda$ in the form $\lambda=\phi(z,\bar z,s)$
 and thus obtain the equation $t =
s^{2b-\ell +1} \phi(z,\bar z,s)$ for $M$.
By the construction of $\phi$, the minimum positive power of $s$
that can appear in its Taylor expansion must be $d-b$.
The desired properties \eqref{desired}
are easily verified from \eqref{E:phi} and \eqref{E:F},
from which it follows that $\Lambda_M$ contains the two
points $(a,a,2b - \ell +1)$ and $(c,a,b - \ell + 1 + d)$, as claimed
above.

Note that for $(\alpha',n',\mu') = (a,a,2b -\ell +1) $ and $(\alpha,n,\mu) = (c,a,b +d -\ell +1)$,
we have
\begin{equation}\label{relation}
\frac{ (\alpha' + n') \mu - (\alpha +n)\mu'}{(\alpha' + n') - (\alpha+n)} = \ell + 1,
\end{equation}
where we have used the relation $c(\ell-b)=a(\ell-d)$.
We claim that,  for $t \in \R$, the biholomorphism
\[ H_t (z,w) = \left( \frac{z}{(1 - t w^{\ell})^{h}},
\frac{w}{(1 - t w^{\ell})^{\frac1\ell}}\right),\]
where $h:=\frac{1}{a}\left( 1- \frac{b}{\ell} \right) = \frac{1}{c}\left( 1- \frac{d}{\ell} \right)>0$,
maps $M$ into itself.
Indeed, it is easy to check that $H_t$ is induced on $M$ by
the biholomorphism
\[ \left( \zeta,\eta \right) \mapsto \left( \frac{\zeta}{1-t\eta}, \frac{\eta}{1-t \eta} \right)
\]
of $S$.
Observe that $j^\ell_0H_t=j^\ell_0 H_{t'}$, for all $t,t'$, but
$j^{\ell+1}_0H_t=j^{\ell+1}_0 H_{t'}$ only if $t=t'$.
Thus the automorphisms of $(M,0)$ are not uniquely determined by their $\ell$-jets.
On the other hand, we have $K=\{\ell+1\}$ by \eqref{relation}
and hence the unique jet determination holds for $(\ell+1)$-jets
in view of Theorem~\ref{T:jetparam}.

\end{document}